\documentclass[10pt,twoside,english]{article}
\usepackage{amssymb,amsmath,babel,geometry,latexsym,graphics,tabularx,shapepar,enumerate}
\usepackage{hyperref}
\usepackage[all,2cell]{xy} \UseAllTwocells \SilentMatrices

\newtheorem{Theorem}{Theorem}
\newtheorem{Lemma}[Theorem]{Lemma}
\newtheorem{Proposition}[Theorem]{Proposition}
\newtheorem{Corollary}[Theorem]{Corollary}
\newtheorem{Definition}[Theorem]{Definition}
\newtheorem{Remark}[Theorem]{Remark}

\newcommand{\ZZ}{\mathbb{Z}}
\newcommand{\FF}{\mathbb{F}}
\newcommand{\CC}{\mathbb{C}}
\newcommand{\QQ}{\mathbb{Q}}
\newcommand{\RR}{\mathbb{R}}

\newcommand{\GG}{\mathbb{G}}

\newcommand{\PP}{\mathbb{P}}

\newcommand\CVD{{\hfill\hfil{\lower 2 pt\hbox{\vrule\vbox to 7pt 
{\hrule width 6pt\vfill\hrule}\vrule}}}\vskip 0.5cm}

\title{On the generalized Carlitz module.
}
\author{Federico Pellarin
}

\begin{document}

\maketitle 

\begin{small}
\noindent\textbf{Abstract.} The aim of this note is to gather formal similarities between two apparently different functions; {\em Euler's  function} $\Gamma$
and {\em Anderson-Thakur function} $\omega$. We discuss these similarities in the framework of the {\em generalized Carlitz's module}, a 
common structure which can be helpful in framing the theories of both functions. We further analyze several noticeable 
differences while investigating the {\em exponential functions} associated to these structures.  
\end{small}

\medskip

\section{Introduction} 

Euler's gamma function 
\begin{equation}\label{gammamellin}\Gamma(s)=\int_0^\infty t^{s-1}e^{-t}dt\end{equation}
is known to satisfy several functional relations, among which the so-called {\em translation formula}, {\em multiplication formulas} and 
the {\em reflection formula}. In a recent work \cite{Pe}, the author used the so-called {\em Anderson-Thakur} function $\omega$ (introduced in \cite{AT})
to deduce properties of a new class of {\em deformations of Carlitz's zeta values}, which are special 
values of zeta functions in the 
framework of the theory of global fields of positive characteristic. Already in \cite{Pe}, we highlighted several similarities 
between $\Gamma$ and $\omega$; in particular, it was observed there that functions closely related to 
$\omega$ also satisfy functional relations similar to those satisfied by $\Gamma$. 

The aim of this expository note is to gather further analogies and highlight differences between these functions
by using the {\em formalism of Carlitz's module} in application to general {\em difference fields} $(\mathcal{K},\tau)$, that is,
fields with a distinguished endomorphism $\tau$.
The idea of making use of the formalism of difference fields is natural and not new, see for example Mumford's paper 
\cite{Mum}. Later, Hellegouarch and his school (see for example \cite{Hell1,Hell2}), highlighted how several 
properties of Carlitz module   (and more generally of Drinfeld modules) are formal consequences of this approach.

A key tool of this formalism is the notion of {\em exponential}. In general, the latter is a formal series of the skew ring of formal series $\mathcal{K}[[\tau]]$.
In many cases, we can associate a function to it. We will discuss the exponential functions underlying 
respectively the gamma function and the
function of Anderson and Thakur (commonly known as the {\em Carlitz's exponential}) and we will also describe how, 
although partially submitted to a common formalism, these 
functions are essentially of different nature.

The paper will also make use of Carlitz's formalism to interpret classical properties of 
the Hurwitz's zeta function; this will be helpful in a tentative to better understand properties of the above mentioned deformations of
Carlitz's zeta values, introduced in \cite{Pe}, and denoted by $L(\chi_t^\beta,\alpha)$. 

\section{The generalized Carlitz's module}


A {\em difference field} $(\mathcal{K},\tau)$ is a field $\mathcal{K}$ with a distinguished field endomorphism $\tau:\mathcal{K}\rightarrow\mathcal{K}$.
A difference field may also be simply denoted by $\mathcal{K}$. 
An extension of difference fields $i:\mathcal{K}\rightarrow\mathcal{K'}$ (alternative notation $\mathcal{K}'/\mathcal{K}$) is an extension of fields equipped with respective endomorphisms $\tau,\tau'$ 
such that $i\tau=\tau'i$ (see Levin's book \cite[Definition 2.1.4]{Levin} for a slightly more general notion and for the background on difference fields). Every difference field can be embedded in an {\em inversive difference field}, that is, a difference field 
such that the distinguished endomorphism is an automorphism (see Cohn's paper \cite{Cohn}; see also Levin, \cite[Proposition 2.1.7]{Levin}). In this case, for all $n\in\ZZ\setminus\{0\}$,
$(\mathcal{K},\tau^n)$ again is an inversive difference field.

The {\em field of constants} $\mathcal{K}^\tau$ of a difference field $\mathcal{K}$ is the subfield of elements $c$ of $\mathcal{K}$
such that $\tau( c)=c$, with $\tau$ the distinguished endomorphism. 
The {\em field of periodic points} $\mathcal{K}^{\text{per.}}$  of $\mathcal{K}$ is the subfield of $\mathcal{K}$
whose elements are the $c\in\mathcal{K}$ such that, for some $n\in\ZZ_{>0}$, we have $\tau^n( c)=c$; this subfield is the algebraic closure of 
$\mathcal{K}^\tau$ in $\mathcal{K}$. 

From now on, all the difference fields we consider are inversive.
The most relevant examples we are interested in, concern the couples $(\mathcal{K},\tau)$ where $\mathcal{K}$ is:
\begin{description}
\item[(1)] The field  $F$ of meromorphic functions $f:\CC\rightarrow\PP_1(\CC)$. In this case, the automorphism 
$$\tau:F\rightarrow F$$ is defined, for $f\in F$ (\footnote{Here, the variable on which $\tau$ acts, will be most of the time denoted by $s$. However, sometimes, 
we will need to deal with several variables and will then adopt the following convention. For a variable $s$, the 
{\em right shift operator} that sends $s$ to $s+1$ and keeps all other variables parameters etc. fixed, will be denoted by $\tau_s$.
For example, we have, for a suitable function of two variables $f(z,s)$, $\tau_s(f(z,s))=f(z,s+1)$.}),
by $$\tau( f)(s)=f(s+1).$$ 
The couple $(F,\tau)$ is an inversive difference field. We have that
$F^\tau$ is the subfield of periodic meromorphic functions and that $F^{\text{per.}}$ is the subfield of 
meromorphic functions which are periodic of integral period, which, as we already pointed out, is algebraically closed in $F$.
\item[(2)] The field $\CC_\infty$ defined as follows. 
Let $p$ be a prime number, $q$ a positive power of $p$ and $\FF_q$ the finite field with $q$ elements. Then $\CC_\infty$
is the completion of an algebraic 
closure of the completion $K_\infty$ of $K=\FF_q(\theta)$ for the $\theta^{-1}$-adic valuation. 
There exists only one absolute value $|\cdot|$ associated to this valuation, such that $|\theta|=q$.
The field $\CC_\infty$
is endowed with the automorphism of infinite order $\tau:\CC_\infty\rightarrow\CC_\infty$ defined, for $c\in\CC_\infty$,
by $$\tau( c)=c^q.$$ The couple $(\CC_\infty,\tau)$ is an inversive difference field. We have $\CC_\infty^\tau=\FF_q$
and $\CC_\infty^{\text{per.}}=\FF_q^{\text{alg.}}$, the algebraic closure of $\FF_q$ in $\CC_\infty$.
\end{description}
In each example, we will try to adopt proper coherent notations. However, this is not always compatible 
with the need for keeping them simple. We will then accept abuses of notation although the reader will be able to
recognize the appropriate structures thanks to a detailed description of the context.

We will mainly focus on these two examples to keep the length of this note reasonable, but many other choices 
of $(\mathcal{K},\tau)$ are interesting; some of them will be occasionally mentioned along this note.

\subsection{Definition and first properties}

We now closely follow Carlitz \cite{Ca0}. Our discussion also follows quite closely the presentation
in Goss' book \cite{Go} (see, for example, Proposition 3.3.10). Here, we additionally deal with a general inversive 
difference field satisfying certain conditions.

Let us fix an inversive difference field $(\mathcal{K},\tau)$ as above. 
We will use the skew ring $\mathcal{K}[[\tau]]$ of formal series 
\begin{equation}\label{series}
\sum_{i\geq 0}c_i\tau^i,\quad c_i\in\mathcal{K},\end{equation}
with the product defined by the formula
$$\sum_{i\geq 0}c_i\tau^i\cdot\sum_{j\geq 0}d_j\tau^j=\sum_{k\geq 0}\left(\sum_{i+j=k}c_i\tau^i(d_j)\right)\tau^k.$$
The ring $\mathcal{K}[[\tau]]$ contains the subring $\mathcal{K}[\tau]$ of 
{\em polynomials in $\tau$}, that is, series (\ref{series}) such that $c_i=0$ for all but finitely many indexes $i$. 
The following construction will look at first sight quite abstract. Eventually, our approach will later become very concrete.

In order to construct generalized Carlitz's modules associated to $\mathcal{K}$ we need to fulfill the following hypothesis: 

\medskip

\noindent\emph{Hypothesis.}
We will suppose, in all the following, that $\mathcal{K}\setminus\mathcal{K}^{\text{per.}}\neq\emptyset$, or in other words, that
$\mathcal{K}$ is a transcendental extension of $\mathcal{K}^{\text{per.}}$,
property true for both the Examples (1) and (2) above. This implies that $\tau$ is of infinite order.

\medskip
 
Let us choose, as we can, $\vartheta\in\mathcal{K}$ transcendental over $\mathcal{K}^{\text{per.}}$ and
let us define, inductively,
for $z\in\mathcal{K}$:
\begin{eqnarray*}
E_0^{(\vartheta)}(z)&=&z\\
E_1^{(\vartheta)}(z)&=&\frac{\tau(E_0^{(\vartheta)}(z))-E_0^{(\vartheta)}(z)}{\tau(\vartheta)-\vartheta}=\frac{\tau(z)-z}{\tau(\vartheta)-\vartheta}\\
E_2^{(\vartheta)}(z)&=&\frac{\tau(E_1^{(\vartheta)}(z))- E_1^{(\vartheta)}(z)}{\tau^2(\vartheta)-\vartheta}\\
& \vdots &\\
E_k^{(\vartheta)}(z)&=&\frac{\tau (E_{k-1}^{(\vartheta)}(z))- E_{k-1}^{(\vartheta)}(z)}{\tau^k(\vartheta)-\vartheta}\\
& \vdots &
\end{eqnarray*} To avoid heavy notations, we will drop the superscript $(\cdot)^{(\vartheta)}$. This construction
is borrowed from Carlitz \cite{Ca0}, see also \cite[Section 3.5]{Go} and \cite[Section 3]{Hell2}.

It is easy to see that the sequence of operators $(E_k)_{k\geq 0}$ is a {\em $\tau$-linear higher derivation} on $\mathcal{K}$.
In other words, for all $k\geq 0$, $E_k$ is $\mathcal{K}^\tau$-linear and if $x,y$ are elements of $\mathcal{K}$,
the following $\tau$-twisted version of Leibniz rule holds, for $k\geq 0$:
\begin{equation}\label{leibniz}
E_k(xy)=\sum_{i+j=k}E_i(x)\tau^{i}(E_j(y)).
\end{equation}

The next definition is then meaningful:
\begin{Definition}\label{definitioncarlitz}
{\em The {\em generalized Carlitz module associated to $(\mathcal{K},\tau)$ and $\vartheta$} is the injective $\mathcal{K}^\tau$-algebra homomorphism
$$\phi_{\mathcal{K},\tau,\vartheta}:\mathcal{K}\rightarrow\mathcal{K}[[\tau]]$$ defined, for $z\in\mathcal{K}$, by
$$\phi_{\mathcal{K},\tau,\vartheta}(z)=\sum_{k\geq 0}(-1)^kE_k^{(\vartheta)}(z)\tau^k.$$}
\end{Definition}

We notice that, independently on the choice of $(\mathcal{K},\tau)$ and $\vartheta$,
$$\phi_{\mathcal{K},\tau,\vartheta}(\vartheta)=\vartheta-\tau.$$ From now on, we will write $\phi$ instead of 
$\phi_{\mathcal{K},\tau,\vartheta}$ to simplify notations, if there is no risk of confusion.
\begin{Lemma}\label{wronskian} For $n\geq 2$,
let $z_1,\ldots,z_n$ be elements of $\mathcal{K}$. We have:
$$
\det\left(\begin{array}{ccc} E_0(z_1) & \ldots & E_0(z_n)\\
\vdots & &\vdots \\
E_{n-1}(z_1) & \ldots & E_{n-1}(z_n)\end{array}\right)=0
$$
If and only if $z_1,\ldots,z_n$ are $\mathcal{K}^\tau$-linearly dependent.
\end{Lemma}
\noindent\emph{Proof.} By induction on $n\geq 2$, one proves the identity
\begin{equation}\label{idFn}
\det\left(\begin{array}{ccc} E_0(z_1) & \ldots & E_0(z_n)\\
\vdots & &\vdots \\
E_{n-1}(z_1) & \ldots & E_{n-1}(z_n)\end{array}\right)=F_n^{-1}\det\left(\begin{array}{ccc} z_1 & \ldots & z_n\\
\vdots & &\vdots \\
\tau^{n-1}(z_1) & \ldots & \tau^{n-1}(z_n)\end{array}\right),\end{equation}
where 
$$F_n=\prod_{i=0}^{n-1}\prod_{j=1}^{n-1-i}(\tau^{i+j}\vartheta-\tau^i\vartheta).$$
Indeed, we can inductively replace the definition of $E_{k}(z_i)$ ($i=1,\ldots,n$) in the following way,
where we have set $G_i=\tau^i\prod_{j=1}^{n-1-i}(\tau^j\vartheta-\vartheta)$ for $i=0,\ldots,n-1$ (so that $F_n=\prod_{i=0}^{n-1}G_i$):
\begin{eqnarray*}
\det\left(\begin{array}{ccc}E_0(z_1)&\ldots&E_0(z_n)\\
E_1(z_1)&\ldots&E_1(z_n)\\
\vdots &  & \vdots \\
E_{n-1}(z_1)&\ldots&E_{n-1}(z_n)\end{array}\right)&=&
\det\left(\begin{array}{ccc}z_1&\ldots&z_n\\
\frac{\tau z_1-z_1}{\tau \vartheta-\vartheta}&\ldots&\frac{\tau z_n-z_n}{\tau \vartheta-\vartheta}\\
\vdots &  & \vdots\\
\frac{\tau E_{n-2}(z_1)-E_{n-2}(z_1)}{\tau^{n-1}\vartheta-\vartheta}&\ldots&\frac{\tau E_{n-2}(z_n)-E_{n-2}(z_n)}{\tau^{n-1}\vartheta-\vartheta}\end{array}\right)\\
&=&G_0^{-1}\det\left(\begin{array}{ccc}z_1&\cdots &z_n\\
\tau z_1 & \cdots & \tau z_n\\ \tau E_1(z_1)&\cdots & \tau E_1(z_n)\\ \vdots & & \vdots \\ \tau E_{n-2}(z_1)&\cdots&\tau E_{n-2}(z_n)\end{array}\right)\\
&=&G_0^{-1}G_1^{-1}\det\left(\begin{array}{ccc}z_1&\cdots &z_n\\
\tau z_1 & \cdots & \tau z_n\\ \tau^2 z_1&\cdots & \tau^2 z_n\\ \tau^2E_1(z_1)& \cdots & \tau^2E_1(z_n)\\ \vdots & & \vdots \\ \tau^2 E_{n-3}(z_1)&\cdots& \tau^2E_{n-3}(z_n)\end{array}\right)=\cdots
\end{eqnarray*}
and so on, hence yielding the identity (\ref{idFn}).

Since $F_n$ is non-zero by the assumption that $\vartheta\not\in\mathcal{K}^{\text{per.}}$, the Lemma follows from (\ref{idFn}) and the Wronskian criterion
(\footnote{The determinants above can also be viewed as generalizations of so-called {\em Moore's determinants} (see \cite[Section 1.3]{Go}); sometimes, they are called {\em casoratians}.}). See also \cite[Lemma 1.3.3]{Go}.\CVD

In all the following, we denote by $\mathcal{A}$ the ring $\mathcal{K}^\tau[\vartheta]\subset\mathcal{K}$. Since 
$$\phi(\vartheta)=\vartheta-\tau,$$
we observe that $\phi$ induces a $\mathcal{K}^\tau$-algebra homomorphism
$$\phi:\mathcal{A}\rightarrow\mathcal{A}[\tau].$$ 
For example,
\begin{eqnarray*}
\phi(\vartheta^2)&=&\vartheta^2-(\vartheta+\tau(\vartheta))\tau+\tau^2,\\
\phi(\vartheta^3)&=&\vartheta^3-(\tau(\vartheta^2)+\vartheta\tau(\vartheta)+\vartheta^2)\tau+(\tau^2(\vartheta)+\tau(\vartheta)+\vartheta)\tau^2-\tau^3.
\end{eqnarray*}
In fact, more is true, as the following elementary proposition shows.

\begin{Proposition}\label{prop1} Let $a$ be an element of $\mathcal{K}$, so that $\phi(a)\in\mathcal{K}[[\tau]]$.
We have that $\phi(a)\in \mathcal{K}[\tau]$ if and only if
$a\in \mathcal{A}$.
\end{Proposition}
\noindent\emph{Proof.} Assume first that $a\in\mathcal{A}=\mathcal{K}^\tau[\vartheta]$. Then, we may write
$$a=\sum_{i=0}^da_i\vartheta^i,\quad a_i\in\mathcal{K}^\tau,\quad a_d\neq0.$$ We then deduce that
$$\phi(a)=a+\cdots+(-1)^da_d\tau^d,$$ agreeing with our statement.

Let us now assume that $\phi(a)\in\mathcal{K}[\tau].$ Then, there exists $d'$ such that
$$\phi(a),\phi(1),\phi(\vartheta),\ldots,\phi(\vartheta^{d'})$$
are $\mathcal{K}$-linearly dependent operators of $\mathcal{K}[\tau]$. By Lemma \ref{wronskian}, $a,1,\vartheta,\ldots,\vartheta^{d'}$ are 
linearly dependent over $\mathcal{K}^\tau$ so that $a\in\mathcal{K}^\tau[\vartheta]$.\CVD

We deduce:

\begin{Corollary} For any choice of $\vartheta\in\mathcal{K}\setminus\mathcal{K}^{\text{per.}}$,
the Carlitz's module induces an injective $\mathcal{K}^\tau$-algebra homomorphism 
$$\phi:\mathcal{A}\rightarrow\mathcal{A}[\tau]$$ such that, for all $a\in\mathcal{A},$ monic of degree $d$ in $\vartheta$, $\phi(a)=a+\cdots+(-1)^d\tau^d$.\end{Corollary}


\medskip

\noindent\emph{Link with Example (1)}. We choose $\vartheta=\text{Id}$, the identity map $\CC\rightarrow\CC$.
This is licit: since $F^{\text{per.}}$ is algebraically closed in $F$ and $\text{Id}$ is not periodic, it is transcendental over $F^{\text{per.}}$.
In the following, we will identify $\text{Id}$ with the variable $s$
of our meromorphic functions in $F$ (yet another abuse of notation). Therefore, we will write $\vartheta=s$.

We recall that for $a\in F$,
$\phi(a)\in F[\tau]$ if and only if $a\in F^\tau[s]$. In particular:
$$\phi(s)=s\tau^0-\tau^1$$
and, for $a\in\CC[s]$ of degree $d$,
$$\phi(a)=\sum_{i=0}^dE_i(a)\tau^i,$$
where $E_i(a)$ is a polynomial of degree $d-i$ in $s$ with complex coefficients. 
For example, one has:
\begin{eqnarray*}
\phi (s^2) & = & s^2 - (2 s + 1)\tau + \tau^2 \\
\phi (s^3) & = & s^3 - (3 s^2 + 3 s + 1)\tau + 
  3 (1 + s)\tau^2 - \tau^3.
\end{eqnarray*}
More generally, one checks easily that if $a$ is of degree $d$ and $c$ is its leading coefficient, 
then $E_d(a)=(-1)^dc$. Our definition allows to evaluate $\phi$ at elements of $\mathcal{K}\setminus\mathcal{A}$ as well. For example,
it is easy to show that:
$$\phi(s^{-1})=\sum_{j\geq 0}(-1)^{j+1}(s)_{j+1}^{-1}\tau^j,$$
where $(s)_j=s(s-1)\cdots(s-j+1)$ denotes the $j$-th Pochhammer polynomial in $s$.

\medskip

\noindent\emph{Link with Example (2)}. If $\vartheta=\theta$,  we shall write, in all the following, $A=\mathcal{A}=\FF_q[\theta]$.
Let us suppose that $\vartheta=\theta$, by Lemma \ref{wronskian}, that $\phi(a)\in \CC_\infty[\tau]$ if and only if 
$a\in A$. More precisely, we have:
\begin{equation}\label{oldcarlitz}
\phi(a)=\sum_{i=0}^d (-1)^iE_i^{(\theta)}(a)\tau^i,\end{equation}
where $d=\deg_\theta(a)$. In particular, $$\phi(\theta)=\theta-\tau,$$ and we obtain
the {\em old} module originally considered by Carlitz, for example in  his papers \cite{Ca0,Ca1,Ca2,Ca3}.

\medskip

\noindent{\bf Important note.} However, to ease the connection with the more recent literature, we will use, from now on in any reference to 
the Example (2), the {\em actual} Carlitz's module, that is,
the $\FF_q$-algebra homomorphism 
$$\phi:\CC_\infty\rightarrow\CC_\infty[[\tau]]$$
determined, for $z\in\mathcal{K}$, by
\begin{equation}\label{newcarlitz}
\phi(z)=\sum_{k\geq 0}E_k^{(\vartheta)}(z)\tau^k.\end{equation}
In particular for $\vartheta=\theta$, 
$$\phi(\theta)=\theta+\tau.$$
This makes our exposition slightly more complicated but should not lead to confusion. In all references to the Example (2) (including when $\vartheta\neq\theta$), Carlitz's 
module 
will be defined by (\ref{newcarlitz}) instead of (\ref{oldcarlitz}). It would be also possible to switch all our discussion 
by considering (\ref{newcarlitz}) to define Carlitz's module for {\em all} inversive difference fields $\mathcal{K}$ as above, but, in this case, we would
lose the direct connection with the theory of the gamma function in the framework of Example (1), as the reader will notice
reading the rest of this paper. Also, we are convinced that Carlitz was aware of these analogies in choosing  (\ref{oldcarlitz}) for his definition and this led us
to make our choice for the present paper.

\subsection{Torsion}

\begin{Definition}{\em Let $a$ be an element of $\mathcal{A}=\mathcal{K}^\tau[\vartheta]$ and $x$ an 
element of $\mathcal{K}$. We say that $x\in\mathcal{K}$ is an element of {\em $a$-torsion} if $\phi(a)(x)=0$.}\end{Definition}
The subset $\text{ker}(\phi(a))$
of elements of $a$-torsion of $\mathcal{K}$ is easily seen to be a $\mathcal{K}^\tau$-vector space 
of finite dimension $\leq d=\deg_\vartheta(a)$ and, at once, an $\mathcal{A}$-module, because if $\phi(a)(x)=0$ then, for all $b\in \mathcal{A}$,
$0=\phi(b)\phi(a)(x)=\phi(a)(\phi(b)(x))$. 

In the framework of Example (1), we have the following result.
\begin{Proposition}\label{Praagman}
Let $f$ be an element of $\mathcal{A}=F^\tau[s]$. The set $\text{ker}(f)$ is a $F^\tau$-vector subspace of $F$
of dimension $d$, where $d$ is the degree in $s$ of $f$.
\end{Proposition} 
\noindent\emph{Proof.} This follows from a more general result due to Praagman, \cite[Theorem 1, p. 102]{Praagman}
which says that, with $\tau$ as in Example (1), a linear homogeneous $\tau$-difference equation of order $d$ with coefficients in $F$ 
can be fully solved in $F$, and the set of solutions is an $F^\tau$-vector subspace of $F$ of dimension $d$.\CVD

In the framework of Example (2), we have the following well known analogue of Proposition \ref{Praagman}.
\begin{Proposition}
Let $a$ be an element of $A=\FF_q[\theta]$ of degree $d$ in $\theta$. The set $\text{ker}(a)$ is an $\mathbb{F}_q$-vector subspace of 
$K^{\text{sep}}$ of dimension $d$.
\end{Proposition} 
We recall that $K^{\text{sep}}$ denotes the maximal separable extension of $K=\FF_q(\theta)$ in $\CC_\infty$.

\medskip

\noindent\emph{Proof.} Let $a$ be as in the statement of the Proposition. Solving the $\tau$-linear equation
$\phi(a)X=0$ amounts to solve a separable algebraic equation of degree $q^d$ the set of solutions of which is 
an $\FF_q$-vector space of dimension $d$.\CVD

In general,
for a given $a\in\mathcal{A}$, the existence of non-zero torsion elements is not guaranteed and reflects the choice 
of the field $\mathcal{K}$. The dimension of the vector space $\text{ker}(a)$ may be, in certain circumstances, strictly smaller that
$d$. But up to extension the difference field $\mathcal{K}$, we may assume that the dimension is exactly $d$.

Indeed, there is a procedure of construction of an extension of $\mathcal{K}$ in which all the spaces $\text{ker}(a)$ have maximal dimension, 
which is similar in spirit to the construction of an algebraically closed extension of a field.
Generalizing algebraic equations in one indeterminate, a $\tau$-{\em difference equation defined over $\mathcal{K}$} is a formal identity 
\begin{equation}\label{differenceequation}
P(X,\tau X,\ldots,\tau^l X)=0,\end{equation}
where $P$ is a non-zero polynomial in $n$ indeterminates. Solving the above $\tau$-difference equation 
amounts to look for the $(l+1)$-tuples 
$$(x_{0},x_{1},\ldots,x_{l})$$ with coordinates in some difference field extension of $\mathcal{K}$ such that
$P(x_0,\ldots,x_l)=0$, and such that $x_{i}=\tau^{i}x_{0}$ for all $i$. Any $x$ with the above properties is called a {\em solution}
of the $\tau$-difference equation (\ref{differenceequation}). 

By the so-called ACFA theory of Chatzidakis-Hrushowski \cite{ChaHru}, there exists an {\em existentially closed}
difference field extension $(\mathbb{K},\tau)$ of $(\mathcal{K},\tau)$. 
We recall here that an existentially closed difference field $(\mathbb{K},\tau)$ 
is characterized by the property that
every finite system of $\tau$-difference equations defined over $\mathbb{K}$ having a solution in some field extension of $\mathbb{K}$
already has a solution in this field. An existentially closed extension $\mathbb{K}$ of $\mathcal{K}$ is, in particular, algebraically closed, 
and such that, for all 
$n>0$, $$\mathbb{K}^{\tau^n}\cap\mathcal{K}=\mathcal{K}^{\tau^n},$$ so that $\mathbb{K}^{\text{per.}}\cap\mathcal{K}=\mathcal{K}^{\text{per.}}$.
There is no canonical choice of such an extension;
this is due to the fact that in general, given two difference field extensions $\mathcal{L}$ and $\mathcal{M}$ of a given 
difference field $(\mathcal{K},\tau)$, a good notion of ``compositum" (a minimal difference field extension of $\mathcal{K}$
in which $\mathcal{L}$ and $\mathcal{M}$ can be simultaneously embedded) is not always available. This yields the
notion of {\em compatible} and {\em non-compatible} difference field extensions. 
See \cite[Chapter 5]{Levin}.

The following Proposition holds.
\begin{Proposition}\label{ssol} Let us assume that $(\mathbb{K},\tau)$ is existentially closed.
If $L=a_0\tau^0+\cdots+a_d\tau^d\in \mathbb{K}[\tau]$ is such that $a_d a_0\neq0$,
then there exist $x_0,\ldots,x_{d-1}\in \mathbb{K}$, $\mathbb{K}^\tau$-linearly independent, such that
$L(x_i)=0$ for $i=0,\ldots,d-1$. Moreover, the set of solutions of the linear $\tau$-difference equation
$L(X)=0$ is equal to the $\mathbb{K}^\tau$-vector space of dimension $d$ generated by $x_0,\ldots,x_{d-1}$.
\end{Proposition}
\noindent\emph{Sketch of proof.} This result is well known, see \cite[Proposition 8.2.4]{Levin}. It is easy to show that any non-trivial linear $\tau$-difference equation 
with coefficients in $\mathbb{K}$ has a non-trivial solution in $\mathbb{K}$.
We proceed by induction on $d\geq 0$. If $d=0$,
the statement is trivial. Let us assume now that $d>0$.
Since $\mathbb{K}$ is existentially closed, there exists a solution $x_0\neq 0$ of $L(x)=0$. Right division algorithm
holds in $\mathbb{K}[\tau]$, so that there exists $\tilde{L}\in \mathbb{K}[\tau]$ unique, with $L=\tilde{L}L_{x_0}$, where, for $x\in \mathbb{K}^\times$, we 
have written $L_x=\tau-(\tau x)/x$. Since the order of the difference operator $\tilde{L}$ is $d-1$, there exist $y_1,\ldots,y_{d-1}$, $\mathbb{K}^\tau$-linearly independent elements of
$\mathbb{K}$ such that $\tilde{L}y_i=0$ for all $i$. Now, for all $i\geq 1$, let $x_i$ be a solution of $L_{x_0}(x_i)=y_i$
(they exist, again because $\mathbb{K}$ is existentially closed). Then, $x_0,x_1,\ldots,x_{d-1}$ are $d$ linearly independent 
elements of $\mathbb{K}$, solutions of $L(x)=0$ such that the set of the solutions in $\mathbb{K}$ of the equation $L(x)=0$,
a $\mathbb{K}$-vector space, has dimension $\geq d$. Now, it is easy to verify, by using the Wronskian Lemma, that the dimension is 
exactly $d$.\CVD
We observe that for all $a\in \mathcal{A}$ of degree $d$ in $\vartheta$, $\phi(a)$ is precisely of the form $L$ as in Proposition \ref{ssol}. 
Therefore, by using Proposition \ref{prop1} and a little computation, we deduce the following Lemma.
\begin{Lemma} Let $a\in\mathcal{A}$ be of degree $d$ in $\vartheta$.
If $\mathcal{K}$ is existentially closed, the vector space $\text{ker}(\phi(a))$ has dimension $d$. 
Moreover, there is an isomorphism of $\mathcal{A}$-modules $\text{ker}(a)\cong \mathcal{A}/a\mathcal{A}$.
\end{Lemma}
We notice that the difference field extension of $\mathcal{K}$ obtained by adjoining the torsion elements is uniquely determined up to difference field isomorphism.

\begin{Definition}
{\em A sequence $(x_i)_{i\geq 1}$ of elements of $\mathcal{K}$ is $\vartheta$-{\em coherent}, or {\em coherent},
if $x_1\neq0$ with $\phi(\vartheta)(x_1)=0$, and if $\phi(\vartheta)(x_i)=x_{i-1}$ for all  $i>1$. Sometimes, for a given coherent sequence 
$\Xi=(x_i)_{i\geq 1}$, we will need to set $x_0=0$.}\end{Definition}
We have the following.
\begin{Lemma}
Let $\Xi=(x_i)_{i\geq 1}$ be a coherent sequence. Then the coefficients $x_1,x_2,\ldots,x_n,\ldots$ are linearly independent over $\mathcal{K}^\tau$.
\end{Lemma}
\noindent\emph{Proof.} Let us suppose by contradiction that for some $a_1,\ldots,a_m\in\mathcal{K}^\tau$ with $a_m\neq0$, we have 
$\sum_{i=1}^ma_ix_i=0$. Then, evaluating $\phi(\vartheta)$ on the left and of the right of this identity and noticing that $\phi(\vartheta)(a_ix_i)=a_i\phi(\vartheta)(x_i)$,
we find $a_2x_1+\cdots+a_mx_{m-1}=0$. More generally, evaluating $\phi(\vartheta^i)$, we obtain at once $$\sum_{i=j}^{m}a_ix_{i-j+1}=0,\quad j=1,\ldots,m.$$
The vector $(x_1,\ldots,x_m)$ is then the unique, trivial solution of a non-singular homogeneous linear system, a contradiction because we have supposed that $x_1\neq0$.\CVD
We consider a new indeterminate $t$ and the field $\mathcal{L}= \mathcal{K}((t))$. We extend $\tau$ to $\mathcal{L}$ by setting
$$\tau\left(\sum_{i\geq i_0}c_it^i\right)=\sum_{i\geq i_0}\tau(c_i)t^i,\quad c_i\in\mathbb{K},\quad i_0\in\ZZ.$$
Then, $\mathcal{L}$ is again an inversive difference field such that, for all $n>0$, $\mathcal{L}^{\tau^n}= \mathcal{K}^{\tau^n}((t))$.
\begin{Definition}{\em 
Let $\Xi=(x_i)_{i\geq 1}$ be a coherent sequence of elements of $\mathcal{K}$. The {\em Akhiezer-Baker series} associated to $\Xi$
is the following element of $\mathcal{K}[[t]]$: $$\omega_{\mathcal{K},\tau,\vartheta,\Xi}(t)=\sum_{i\geq 0}x_{i+1}t^i.$$
To simplify our notations, we will write $\omega_{\Xi}$ instead of $\omega_{\mathcal{K},\tau,\vartheta,\Xi}$ when the reference to
the data $\mathcal{K},\tau,\vartheta$ is clearly indicated.}\end{Definition} 
If $a=a(\vartheta)$ is an element of $\mathcal{A}$, we will write $a(t)$ for the function of the variable $t$ or the element of
$\mathcal{L}$ obtained by formal replacement of $\vartheta$ with $t$, which is meaningful because $\vartheta$ is transcendental 
over $\mathcal{K}^\tau$. The following Proposition justifies the adopted terminology and is easy to prove.

\begin{Proposition}\label{PsiX} Let $\Xi=(x_i)_{i\geq 1}$ be a sequence of $\mathcal{K}$ with $x_1\neq0$ and let us 
write $\omega_\Xi=\sum_{i\geq 0}x_{i+1}t^i$. The following properties are equivalent.
\begin{enumerate}
\item $\Xi$ is coherent,
\item $\phi(\vartheta)(\omega_\Xi)=t\omega_\Xi$,
\item For all $a\in\mathcal{A}$, $\phi(a)\omega_\Xi=a(t)\omega_\Xi$.
\end{enumerate}
In particular, if $\Xi,\Xi'$ are two coherent sequences, then
$$\omega_{\Xi}=\lambda_{\Xi,\Xi'}\omega_{\Xi'},$$ where $\lambda_{\Xi,\Xi'}$ is an
element of $\mathcal{L}^\tau$ determined by $\Xi,\Xi'$.
A coherent sequence can be always found in an appropriate extension of the difference field $\mathcal{K}$.
If $\mathcal{K}$ contains a coherent sequence, then it contains all the coherent sequences.
\end{Proposition}
\noindent\emph{Sketch of proof.} The equivalence of points 1 and 2 follows from the identities
$$\phi(\vartheta)(\omega_\Xi)=\sum_{i\geq 0}\phi(\vartheta)(x_{i+1})t^i,\quad t\omega_\Xi=\sum_{i\geq 1}x_it^i.$$
The point 3 clearly implies the point 2. The opposite implication is obtained observing that if $a=a_0+a_1\vartheta+\cdots+a_d\vartheta^d$
with $a_0,\ldots,a_d\in\mathcal{K}^\tau$, then $\phi(a)\omega_\Xi=a_0\omega_\Xi+a_1\phi(\vartheta)\omega_\Xi+\cdots+a_d\phi(\vartheta^d)\omega_\Xi$,
so that $\phi(a)\omega_\Xi=a(t)\omega_\Xi$. The remaining parts of the Proposition follow from the fact that 
the set of solutions in $\mathcal{L}$ of the difference equation $\phi(\vartheta)X=tX$ is either trivial, either a $\mathcal{L}^\tau$-vector space
of dimension $1$.\CVD

In particular, if $\mathcal{K}$ is existentially closed, it contains all the coherent sequences.
\begin{Remark}{\em 
The link with the classical Akhiezer-Baker functions of Krichever's axiomatic approach \cite{Kri} lies in the fact that,
just as the latter, the functions $\omega_{\mathcal{K},\tau,\vartheta,\Xi}$ are eigenfunctions of the full set
of operators $\phi_{\mathcal{K},\tau,\vartheta}(a)\in\mathcal{K}[[\tau]]$, with $a\in \mathcal{A}$. 
Akhiezer-Baker functions are usually defined over 
Riemann surfaces and are characterized by certain essential singularities. 
Baker, in \cite{Bak}, noticed the possibility to relate them to the theta series associated to the surface,
and Akhiezer pointed out that they can also be constructed, in certain cases, as eigenfunctions of
differential operators of order two.

In order to make the link with the classical 
theory  and our choice of terminology clearer we give here a suggestive example of a very particular family of Akhiezer-Baker functions, eigenfunctions
of a Lam\'e operator. The example we give
is treated by N. I. Akhiezer in his book \cite[Chapter 11]{Akh} (see also Krichever's paper \cite{Kri0}). 

Let $E$ be a complex elliptic curve with Weierstrass' model $y^2=4x^3-g_2x-g_3$, analytically isomorphic to
a complex torus $\CC/\Lambda$, where $\Lambda=\ZZ\omega_1+\ZZ\omega_2$ is its lattice of periods.
Let $\wp$ denote the Weierstrass function associated to $\Lambda$, so that, for $\nu\not\in\Lambda$,
$x=\wp'(\nu)$ and $y=\wp(\nu)$. We also denote by $\zeta$ the Weierstrass zeta function of $\Lambda$ and by $\sigma$ the 
Weierstrass sigma function of $\Lambda$.
Let us choose a point $(x^*,y^*)\in E(\CC)$ with $(x^*,y^*)=(\wp(\nu^*),\wp'(\nu^*))$, for some
$\nu^*\in\CC\setminus\Lambda$ and let us additionally consider a complex parameter $z$. The Akhiezer-Baker function
associated to the datum $(x^*,y^*,z)$ is the function of the variables $(x,y)\in E(\CC)$ depending on the 
parameter $z$ defined by
$$\Omega(x,y,z)=\frac{\sigma(\nu+\nu^*-z)e^{z\zeta(\nu)}}{\sigma(\nu-\nu^*)\sigma(\nu^*+z)}.$$
The reader can check, with the elementary properties of $\zeta$, $\sigma$ in mind, that the above function
is a well defined meromorphic function on $E(\CC)\setminus\{\infty\}$ but is not an elliptic function, the obstruction being an (unique) essential singularity at infinity
(that is, at $\nu\equiv0\pmod{\Lambda}$). By the way, it also has a pole at $(x^*,y^*)$. So far, we have looked at this function
in the dependence of the variables $(x,y)$. Now, one verifies easily,
with the help of the classical differential properties of basic Weierstrass functions, that $\Omega(x,y,z)$, as a function
of the variable $z$, for any choice of $x,y,x^*,y^*$, is eigenfunction of the Lam\'e differential operator
$\frac{\partial^2}{\partial z^2}-2\wp(z)$ with eigenvalue $-\wp(\nu)$, namely:
$$\left(\frac{\partial^2}{\partial z^2}-2\wp(z)\right)\Omega(x,y,z)=-\wp(\nu)\Omega(x,y,z).$$ Therefore, this function is eigenfunction
of all the operators in the image of the Krichever module uniquely determined by the above operator, and the corresponding eigenvalues
are obtained from the constant terms of the operators, with the variable $z$ replaced by the parameter $\nu$ (upon certain necessary normalizations).
In analogy with our semi-character $a(\theta)\mapsto a(t)$, this suggests the terminology adopted in the present paper.}\end{Remark}

\subsection{The gamma function as a torsion element} 

The function $\Gamma$ is traditionally defined as the Mellin transform (\ref{gammamellin}) of the function $e^{-z}$.
Since the relation $\Gamma(s+1)-s\Gamma(s)=0$ holds, this function belongs to the $F^\tau$-vector space of $s$-torsion for the 
generalized Carlitz module $\phi$ associated to the difference field $(F,\tau)$. 
More generally, we have the following proposition.
\begin{Proposition}\label{torsion} The sequence (\footnote{The notation $f^{(n)}(s)$ denotes, all along this paper, 
the derivative $\frac{d^n}{ds^n}f(s)$.}) $$\Xi=\left((-1)^{i-1}\frac{\Gamma^{(i-1)}}{(i-1)!}\right)_{i\geq 1}$$
is coherent.
For all $k\geq 1$, the $F^\tau$-vector space $\text{ker}(\phi(s^k))$ of the solutions $X$ of the linear $\tau$-difference equation
$$\phi(s^k)(X)=0$$
has dimension $k$ and is spanned by the $F^\tau$-linearly independent functions
$$\Gamma,\Gamma',\ldots, \Gamma^{(k-1)}.$$
\end{Proposition}
\noindent\emph{Proof.} Let $f$ be a function of $F$ and $a=a(s)$ and element of $F^\tau[s]$. Then we can formally compute,
with the rule $\frac{d}{ds}\tau=0$ and for $n>0$,
$$\frac{d}{ds}(\phi(a^n)f(s))=n\left(\frac{d}{ds}\phi(a)\right)\phi(a^{n-1})f(s)+\phi(a^{n})\frac{d}{ds}f(s).$$
Since
$$\frac{d}{ds}\phi(s)=1,$$
we have in particular that
$$\frac{d}{ds}(\phi(s^n)f(s))=n\phi(s^{n-1})f(s)+\phi(s^{n})\frac{d}{ds}f(s)$$
and if $f(s)$ is $s^n$-torsion, that is, $\phi(s^n)f(s)=0$, we get
$$\frac{\phi(s^n)f'(s)}{n}=-\phi(s^{n-1})f(s).$$
By induction on $n$ we then see that if $\phi(s)f=0$, then $(-1)^n\frac{f^{(n)}}{n!}\in\text{ker}(\phi(s^{n+1}))$. On the other 
hand, it is easy to verify the following formal identity:
$$\frac{1}{m!}\left(\frac{d}{ds}\right)^m\phi(s^n)=\binom{n}{m}\phi(s^{n-m}),$$
from which one deduces, applying the operator $(d/ds)^n$ on both sides of the identity $\phi(s)f=0$,
that $\phi(s)f^{(n)}=nf^{(n-1)}$. The above discussion holds if for example $f=\Gamma$, which yields the first properties of the sequence $\Xi$.
Notice that the linear independence of the coefficients of $\Xi$ can be also verified directly by comparing the orders of 
the poles of the functions $\Gamma^{(n)}$. In particular, $F$ contains all the coherent sequences for $\phi_{F,\tau,s}$.\CVD

In the present framework, we have, after Proposition \ref{PsiX}:
\begin{Proposition}\label{likeomega} The Akhiezer-Baker function associated to the coherent 
sequence $\Xi$ of Proposition \ref{torsion} is the formal series
$$\omega_{\Xi}(t)=\Gamma(s-t)=\sum_{k\geq 0}(-1)^k\frac{\Gamma^{(k)}(s)}{k!}t^k.$$
\end{Proposition}
One can further verify that the series above is convergent for complex numbers $s,t$ such that $|t|<|s|<1$, $\Re(s)>0$.




\begin{Remark}
\label{digamma}
{\em An interesting property following from Proposition \ref{likeomega} is that the function
$$\psi(s-t)=\frac{\Gamma'(s-t)}{\Gamma(s-t)}=\sum_{n\geq0}(-1)^n\frac{\psi^{(n)}(s)}{n!}t^n,$$ where $\psi^{(n)}(s)=(d/dz)^{n+1}\log\Gamma(s)$ (\footnote{Here $\log$ denotes the 
principal determination of the classical logarithm.}) 
is the $n$-th polygamma function, is  a solution of the difference equation
\begin{equation}\label{taudigamma}
\tau X-X=\frac{1}{s-t}.\end{equation} }\end{Remark}

\begin{Remark}{\em We mention that the torsion of the generalized Carlitz's module 
is related to classical functions also in the framework of other difference fields $(\mathcal{K},\tau)$. The next is an important example
that
the reader can find of parallel interest (we might refer to it as to the {\em Example (3)}; see also \cite[pp. 350-354]{Hell2}). Consider $\mathcal{K}$ the field of meromorphic functions over $\CC^\times$ of which the variable is denoted by $x$, $q$ a non-zero complex number 
such that $|q|>1$ and 
$\tau$ the automorphism of $\mathcal{K}$ defined by $\tau x=qx$, so that $\mathcal{K}^\tau$ is the field of elliptic functions over the 
elliptic curve $\CC^\times/q^\ZZ$. Then, Carlitz's formalism applies to $(\mathcal{K},\tau)$
and it turns out that {\em Jacobi's theta series}:
$$x_1=\sum_{m\in\ZZ}q^{-\frac{m(m+1)}{2}}x^m$$
represents a generator of the $x$-torsion. More generally, one can prove that, by setting
$$d_m^{[n]}=q^{-\frac{1}{2}m(m+2n+1)}\prod_{i=1}^n\frac{q^{i+m}-1}{q^i-1}$$
and
$$x_{n+1}=\sum_{m\in\ZZ}d^{[n]}_mx^m,\quad n\geq 0,$$
the sequence $\Xi=(x_n)_{n\geq 0}$ (with $x_0=0$) is coherent, from which one obtains explicitly the corresponding Akhiezer-Baker function
$\omega_\Xi$.}\end{Remark}

\subsection{Some elements of $\text{ker}(\phi(f))$ for $f\in F^\tau[s]$}

Let $f$ be an element of $F^\tau[s]$. We have seen in Proposition \ref{Praagman} that the set
$\text{ker}(f)$ is a $F^\tau$-vector subspace of $F$ of dimension $\deg_s(f)$. Here, we give some
examples of elements of $\ker(f)$.

By Proposition \ref{likeomega}, for all $f\in F^\tau[s]$, 
$$\phi(f)\Gamma(s-t)=f(t)\Gamma(s-t).$$
If $f\in\CC[s]$, and if $x\in\CC$ is a root of $f$, then, obviously,
$$
\phi(f)\Gamma(s-x)=f(x)\Gamma(s-x)=0,
$$
and the function $\Gamma(s-x)$, meromorphic on $\CC$, lies in the kernel of $\phi(f)$. 
This simple argument can be deduced from the case of $f$ of degree one in $s$. Indeed, $\Gamma(s-x)$
obviously generates the kernel of $\phi(s-x)$ (one uses that the field $\CC$ is
algebraically closed).

\medskip

\noindent\emph{Example.} The functions $\Gamma(s-i)$ and $\Gamma(s+i)$ generate the two-dimensional kernel of $\phi(s^2+1)$.

\medskip

If we consider now $f=s-x(s)$ with $x\in F^\tau$ not necessarily constant,
then, away from a discrete subset of $\CC$ depending on $x$, the function
$$\Gamma(s-x(s))$$ is well defined and holomorphic. By the above observations, it certainly generates the kernel of $\phi(f)$.
This already allows to compute the kernel of $\phi(f)$ when $f=\prod_{i=1}^d(s-x_i(s))$, where $x_i\in F^\tau$ are distinct:
the kernel is generated by the functions 
$$\Gamma(s-x_1(s)),\ldots,\Gamma(s-x_d(s)),$$ which are linearly independent over $F^\tau$.
In fact, by using arguments as in the proof of Proposition \ref{likeomega}, we get the following simple improvement. The details of
the proof are left to the reader.

\begin{Proposition}\label{sometorsion} Let $f=\prod_{i=1}^d(s-x_i(s))^{k_i}$ be a monic polynomial of $F^\tau[s]$ of degree $\sum_{i=1}^dk_i$ such that the 
roots $x_i$ are distinct elements of $F^\tau$. Then, the kernel of $\phi(f)$ is spanned by the $F^\tau$-linearly independent functions
$$\Gamma(s-x_1(s)),\ldots,\Gamma^{(k_1-1)}(s-x_1(s)),\ldots,\Gamma(s-x_d(s)),\ldots,\Gamma^{(k_d-1)}(s-x_d(s)),$$
holomorphic on the complement in $\CC$ of a discrete set depending on $f$.\end{Proposition}

The above proposition can be applied when $f$ splits as a product of linear polynomials but cannot applied in the opposite
situation, when $f$ is irreducible as a polynomial in $s$.
If $f\in F^\tau[s]$ is irreducible, there exists a unique irreducible $g\in F^\tau[X]$
such that $g(s)=f$. Assuming further that the roots $x_i(s)$ of $g$ are in $F^{\text{per.}}$, the functions $\Gamma(s-x_i(s))$ are 
well defined on an open subset of $\CC$ as above, but there is no reason that they are 
in the kernel of $\phi(f)$. 

However, $\tau$ induces a permutation $\sigma$ of the $x_i$. Let $t_1,\ldots,t_d$ be independent 
variables such that $d=\deg f$ and such that $\tau t_i=t_i$ for all $i$.
If we set
$$G=G(s,t_1,\ldots,t_s)=\sum_i\Gamma(s-t_i),$$
we have
\begin{eqnarray*}
(\tau G)|_{t_i\mapsto x_i(s)}&=&\sum_i\Gamma(s+1-x_i(s+1))\\
&=&\sum_i\Gamma(s+1-x_{\sigma(i)}(s)))\\
&=&\tau\left(G|_{t_i\mapsto x_i(s)}\right).
\end{eqnarray*}
Hence,
\begin{eqnarray*}
\phi(f)(\sum_i\Gamma(s-x_i(s)))&=&(\phi(f)G)|_{t_i\mapsto x_i(s)}\\
&=&\left(\sum_i f(t_i)\Gamma(s-t_i)\right)|_{t_i\mapsto x_i(s)}\\
&=&0
\end{eqnarray*}
and we get the following Proposition.
\begin{Proposition}
Let $f\in F^\tau[s]$ be irreducible of degree $d$, such that all its roots $x_i(s)$ are in $F^{\text{per}}$. Then, the function
$$\sum_i\Gamma(s-x_i(s)),$$ defined on the complement in $\CC$ of a discrete subset depending on $f$,
is a non-trivial element of the kernel of $\phi(f)$.
\end{Proposition}

\subsection{The basic functional relations}\label{funcrel}

So far, we have considered separately the Carlitz's modules associated to individual data $(\mathcal{K},\tau,\vartheta)$,
but we may well compare different structures over the same field $\mathcal{K}$. For example, we can
choose Carlitz's module associated to the difference field $(\mathcal{K},\tau^n)$ with $\vartheta$ fixed and $n$ varying, or vary $\vartheta$ while fixing $(\mathcal{K},\tau)$. 
Here, we show that the Akhiezer-Baker functions
corresponding to these choices satisfy functional relations that in the framework of Example (1) are at the origin of the classical functional relations
of Euler's gamma function.

\subsubsection{Multiplication relations}

We compare the structures of the Carlitz modules $\phi_{\mathcal{K},\tau^n,\vartheta}$ and 
$\phi_{\mathcal{K},\tau,\vartheta}$ through their Akhiezer-Baker functions. from now on, we assume that $\mathcal{K}$ is large enough
to contain all the coherent sequences associated to the various Carlitz's modules $\phi_{\mathcal{K},\tau,\vartheta}$ that we are considering;
this happens if, for instance, $\mathcal{K}$ is existentially closed.
Let us notice that if $F$ is a non-zero solution of the difference equation 
\begin{equation}\label{FAF}
\tau F=AF\end{equation} for some $A\in\mathcal{K}^\times$, then
\begin{equation}\label{multiplicationF}
\tau^nF=(\tau^{n-1}A)(\tau^{n-2}A)\cdots(\tau A)A F.\end{equation}
Let us now suppose that there exists  a solution $G\in\mathcal{K}^\times$ of the difference equation
\begin{equation}\label{GAG}
\tau^nG=AG.\end{equation}
Then, for all $i=0,\ldots,n-1$, $\tau^n(\tau^iG)=(\tau^iA)(\tau^iG)$ so that
the product
$$H=G(\tau G)\cdots(\tau^{n-1}G)$$ is again solution of (\ref{multiplicationF}). The ratio of any two non-zero solutions of the equation (\ref{multiplicationF})
is in $\mathcal{K}^{\tau^n}$. Therefore,
\begin{equation}\label{GH}
F=\lambda H,
\end{equation} for some $\lambda\in\mathcal{K}^{\tau^n}$.

The above arguments hold for any choice of $A\in\mathcal{K}^\times$. We now focus on the case $A=\vartheta-t$.
In this case, any solution $F$ of (\ref{FAF}) is a multiple by an element of $\mathcal{K}^\tau[[t]]$ of 
an Akhiezer-Baker function $\omega_{\mathcal{K},\tau,\vartheta,\Xi}$ associated to a coherent sequence $\Xi$.

On the other hand, any solution $G$ of (\ref{GAG}) is a multiple by an element of $\mathcal{K}^{\tau^{n}}[[t]]$
of an Akhiezer-Baker function $\omega_{\mathcal{K},\tau^n,\vartheta,\Xi}$ associated to a coherent sequence $\Xi'$.
Therefore, looking at (\ref{GH}) we obtain the following proposition.
\begin{Proposition}[Multiplication relation for Akhiezer-Baker functions]\label{propmultiplication}
Consider a coherent sequence $\Xi$ for $\phi_{\mathcal{K},\tau,\vartheta}$ and a coherent sequence $\Xi'$ for $\phi_{\mathcal{K},\tau^n,\vartheta}$,
$n$ being a positive integer. If we set $\omega=\omega_{\mathcal{K},\tau,\vartheta,\Xi}$
and $\omega'=\omega_{\mathcal{K},\tau^n,\vartheta,\Xi'}$ for the corresponding Akhiezer-Baker functions, then, there exists a non-zero element $\lambda\in\mathcal{K}^{\tau^n}[[t]]$
such that
$$\omega=\lambda(\tau^{n-1}\omega')(\tau^{n-2}\omega')\cdots(\tau\omega')\omega'.$$
\end{Proposition}

\begin{Remark}{\em The well known 
Gauss' {\em multiplication formulas} for Euler's gamma function
can be deduced by first applying Proposition \ref{propmultiplication} and then computing the 
function $\lambda$ by means of Stirling asymptotic formula 
\begin{equation}\label{stirling}
\Gamma(s)\sim\sqrt{2\pi}s^{s-1/2}e^{-s},\quad |\arg(s)|<\pi.\end{equation} An alternative way to prove this property is to
use the digamma function by computing a logarithmic derivative, observing that
the functions $$G(s)=\sum_{k=0}^{N-1}\psi\left(s+\frac{k}{N}\right)$$ (for fixed $N>0$) and
$N\psi(Ns)$ are both solutions of the $\tau$-difference equation
$$\tau X-X=\sum_{k=0}^{N-1}\frac{N}{Ns+k},$$ so that they differ by an element of $F^\tau$.
The details are left to the reader.
Notice also that the function $\lambda$ of Proposition \ref{propmultiplication} depends on $t$.}\end{Remark}

\subsubsection{Cyclotomic relations}

Let $\zeta\in\mathcal{K}^\tau$ be a root of unity of order $n>0$. We can compare the Akhiezer-Baker 
functions of the data $(\mathcal{K},\tau,\vartheta\zeta^i)$, ($i=0,\ldots,n-1$) and $(\mathcal{K},\tau,\vartheta^n)$
by using the obvious identity:
\begin{equation}\label{cyclotomicrel}
\prod_{i=0}^{n-1}(t-\zeta^i\vartheta)=t^n-\vartheta^n.
\end{equation}
Indeed, let us choose suitable coherent sequences $\Xi_0,\ldots,\Xi_{n-1}$ and $\Xi$ associated respectively to the Carlitz's
modules $\phi_{\mathcal{K},\tau,\vartheta}, \phi_{\mathcal{K},\tau,\zeta\vartheta},\ldots\phi_{\mathcal{K},\tau,\zeta^{n-1}\vartheta}$ and 
$\phi_{\mathcal{K},\tau,\vartheta^n}.$
Then, the corresponding Akhiezer-Baker functions satisfy:
\begin{eqnarray*}
\tau\omega_{\mathcal{K},\tau,\zeta^i\vartheta,\Xi_i}(t)&=&(t-\zeta^i\vartheta)\omega_{\mathcal{K},\tau,\zeta^i\vartheta,\Xi_i}(t),\quad i=0,\ldots,n-1,\\
\tau\omega_{\mathcal{K},\tau,\vartheta^n,\Xi}(t^n)&=&(t^n-\vartheta^n)\omega_{\mathcal{K},\tau,\vartheta^n,\Xi}(t^n).
\end{eqnarray*}
Multiplying term-wise the first $n$ identities, we obtain the following proposition.
\begin{Proposition}[Cyclotomic relations for Akhiezer-Baker functions]\label{propcyclotomic} In the above no\-tatio\-ns, there exists a non-zero
element $\mu\in\mathcal{K}^{\tau}[[t]]$ such that
$$\prod_{i=0}^{n-1}\omega_{\mathcal{K},\tau,\zeta^i\vartheta,\Xi_i}(t)=\mu \omega_{\mathcal{K},\tau,\vartheta^n,\Xi}(t^n).$$
\end{Proposition}

The {\em reflection formula} for Euler's gamma function is not directly affiliated to cyclotomic relations.
The appropriate framework to look at it seems to be that of the {\em adjoint} of Carlitz's module.

\subsubsection{Adjunction relations}

Here, we look at the Carlitz module $\phi_{\mathcal{K},\tau^{-1},\vartheta}$. This is different from the modules 
analyzed above because its image is in the new ring of formal series $\mathcal{K}[[\tau^{-1}]]$, but it can be handled in a similar way. 
Given a coherent 
sequence $\Xi$ for $\phi_{\mathcal{K},\tau,\vartheta}$, and setting $\psi=\frac{1}{\tau\omega_{\mathcal{K},\tau,\vartheta,\Xi}}$,
from 
$$\frac{1}{\tau\omega_{\mathcal{K},\tau,\vartheta,\Xi}}=\frac{1}{(t+\vartheta)\omega_{\mathcal{K},\tau,\vartheta,\Xi}}$$
we deduce 
$$\tau^{-1}\psi=(t+\vartheta)\psi.$$
By Proposition \ref{PsiX}, there exists a coherent sequence  $\Xi'$ for $\phi_{\mathcal{K},\tau^{-1},\vartheta}$, the {\em adjoint sequence}, 
such that
$$\omega_{\mathcal{K},\tau^{-1},\vartheta,\Xi'}=\frac{1}{\tau \omega_{\mathcal{K},\tau,\vartheta,\Xi}}.$$

Additionally, in the special case of Example (1), where the difference field is $(F,\tau)$ and $\vartheta=s$, we can choose
$\Xi$ as in Proposition \ref{torsion} so that $\omega_{\mathcal{K},\tau,\vartheta,\Xi}=\Gamma(s-t)$.
We observe that, if $\tau_t$ is the $F$-linear operator (on a suitable space of functions of the variable $t$) which sends $t$ to $t+1$, then
\begin{equation}\label{particularlysimple}\tau(s-t)=\tau_t^{-1}(s-t).\end{equation}
This means that $\Gamma(t-s)$ is an Akhiezer-Baker function for $\phi_{\mathcal{K},\tau^{-1},\vartheta}$ and 
summing everything together,
$$\Gamma(t-s)=\frac{\nu}{\Gamma(s-t+1)},$$ for some $\nu\in\mathcal{K}^{\tau^2}[[t]]$. Of course 
these arguments are covered by simply verifying that the function $\Gamma(s)\Gamma(1-s)$ is 
periodic of period $2$, which is trivial, but our discussion helps in observing that the crucial point is 
(\ref{particularlysimple}), which does not hold for general difference fields $\mathcal{K}$. In the case of Example (2), this identity is replaced
by $$\tau(\vartheta-t)=\tau_t^{-1}(\vartheta-t)^q,$$ where $\tau_t$ is $\CC_\infty$-linear such that $\tau_tt=t^q$.
We do not know whether this can be used to exhibit some kind of adjunction relation for the Akhiezer-Baker functions
in this case, but multiplication relations and cyclotomic relations in this case hold and will be described in Section \ref{functionalrelations}.

\subsection{Exponential and logarithm}

We introduce now another important tool of the theory, the {\em exponential} and the {\em logarithm} series associated to 
the data $(\mathcal{K},\tau,\vartheta)$. These will be 
elements of $\mathcal{K}[[\tau]]$ and we will see that sometimes and quite naturally, it is possible to associate to these formal series certain operators.

\begin{Proposition}
For any given inversive difference field $(\mathcal{K},\tau)$ and $\vartheta$ as above, there exist, unique,
two series $E=E_{\mathcal{K},\tau,\vartheta},L=L_{\mathcal{K},\tau,\vartheta}\in\mathcal{K}[[\tau]]$ with the following properties.
\begin{enumerate}
\item The series are normalized, that is, $E=\tau^0+\cdots$, $L=\tau^0+\cdots$.
\item The series are one inverse of the other for the product rule of $\mathcal{K}[[\tau]]$: $EL=LE=\tau^0.$
\item For all $z\in\mathcal{K}$, we have
$\phi(z)E=E z.$
\item For all $z\in\mathcal{K}$, we have
$L\phi(z)=zL.$
\end{enumerate}
\end{Proposition}
\noindent\emph{Proof.} It follows from an elementary computation. The coefficients of $E$ can be inductively 
computed from the identity $\phi(\vartheta)E=(\vartheta\tau^0-\tau)E=E\vartheta$:
$$E=\sum_{n\geq 0}d_n^{-1}\tau^n,$$
where $d_0=1$ and
\begin{equation}\label{recursionE}
d_n=(\vartheta-\tau^n\vartheta)\tau d_{n-1}, \quad n\geq 1.
\end{equation}
Similarly, the series $L$
$$L=\sum_{n\geq 0}l_n^{-1}\tau^n$$
is uniquely determined by the recursion
\begin{equation}\label{recursionL}l_n=(\tau^n\vartheta-\vartheta)l_{n-1},\quad n\geq 1,\end{equation} with $l_0=1$. Then,
one notices, for $f\in \mathcal{K}$, the identity
$\phi(f)=EfL$.\CVD

In the case of Example (1), where the difference field is $(F,\tau)$, one immediately verifies that
$$E=\sum_{n\geq 0}(-1)^n\frac{\tau^n}{n!},\quad L=\sum_{n\geq 0}\frac{\tau^n}{n!}.$$ 
More generally, we will need the {\em modified logarithm}. Let $x$ be a complex number. The modified logarithm $L_x$
is the element of $F[[\tau]]$:
$$L_x=\sum_{n\geq 0}\frac{x^n\tau^n}{n!}.$$
Obviously, $E=L_{-1}$, $L=L_1$ and $L_xL_{-x}=L_{-x}L_x=\tau^0$. Another useful identity, holding in $F[[\tau]]$, is
$$L_xs=\phi_x(s)L_x,$$ where $\phi_x(s)=1+x\tau$.

\section{Link with Example (2): Anderson-Thakur's functions}

In the case of Example (2), where by convention we have chosen (\ref{newcarlitz}) to define Carlitz's module
(in particular, $\phi_{\CC_\infty,\tau,\vartheta}(\vartheta)=\vartheta+\tau$), the coefficients $d_n,l_n$ appearing in the expansions of 
the series $E_{\CC_\infty,\tau,\vartheta},L_{\CC_\infty,\tau,\vartheta}$ are polynomials of $\vartheta$ and are given by the following formulas:
\begin{eqnarray}
d_n(\vartheta)&=&(\vartheta^{q^n}-\vartheta^{q^{n-1}})(\vartheta^{q^n}-\vartheta^{q^{n-2}})\cdots(\vartheta^{q^n}-\vartheta),\nonumber\\
l_n(\vartheta)&=&(-1)^n(\vartheta^{q^n}-\vartheta)(\vartheta^{q^{n-1}}-\vartheta)\cdots(\vartheta^{q}-\vartheta).\label{dis}
\end{eqnarray}
From now on, we shall write, when dealing with Example (2),
$$\exp_{\tau,\vartheta}=E_{\CC_\infty,\tau,\vartheta},\quad \log_{\tau,\vartheta}=L_{\CC_\infty,\tau,\vartheta}.$$
When $\vartheta=\theta$, we will write $\exp=\exp_{\tau,\theta}$, $\log=\log_{\tau,\theta}$. 
The properties of the latter functions are described in detail in Goss' book \cite{Go} and 
in Thakur's book \cite{Tha}. But when $|\vartheta|>1$, the arguments of these authors can be easily 
transposed to the more general setting of the functions $\exp_{\tau,\vartheta}$, $\log_{\tau,\vartheta}$.
This follows from the elementary observation that we then have, in (\ref{dis}), $|d_n(\vartheta)|=|\vartheta|^{nq^n}$ and
$|l_n(\vartheta)|=|\vartheta|^{\frac{q(q^n-1)}{q-1}}$ for all $n\geq 0$ (see \cite[Sections 3.2, 3.4]{Go}).
So there is no surprise in finding the next Theorem, which specializes in well known results if $\vartheta=\theta$.
The original description involving the {\em old} Carlitz module when $\vartheta=\theta$ can be found in \cite{Ca0}.

\begin{Theorem}\label{expcarlitz} Assuming that $\vartheta\in\CC_\infty$ is such that $|\vartheta|>1$,
The series $\exp_{\tau,\vartheta}\in \FF_q(\vartheta)[[\tau]]$ defines an entire surjective $\FF_q$-linear function, again denoted by
$\exp_{\tau,\vartheta}$:
$$\exp_{\tau,\vartheta}:\CC_\infty\rightarrow\CC_\infty$$ by $\exp_{\tau,\vartheta}(z)=\sum_{i\geq 0}d_i(\vartheta)^{-1}z^{q^i}$. The kernel of this function
is a discrete, free $\mathcal{A}$-submodule of $\CC_\infty$ of rank one.
\end{Theorem}
We choose, once and for all, a generator $\widetilde{\pi}_{\tau,\vartheta}$ of the kernel of $\exp_{\tau,\vartheta}$. It can be proved (see \cite{Go} for the case $\vartheta=\theta$, the general case can be proved by following the same methods) that,
up to a choice of a $(q-1)$-th root of $-\vartheta$,
\begin{equation}\label{tildepidef}
\widetilde{\pi}_{\tau,\vartheta}=\vartheta(-\vartheta)^{\frac{1}{q-1}}\prod_{i=1}^\infty(1-\vartheta^{1-q^i})^{-1}\in (-\vartheta)^{\frac{1}{q-1}}\FF_q((\vartheta^{-1})).
\end{equation}
If $\vartheta=\theta$, we shall write $\widetilde{\pi}$ instead of $\widetilde{\pi}_{\tau,\vartheta}$. This is the so-called 
{\em fundamental period} of Carlitz's module in the framework of Example (2) (there are $q-1$ such choices). It is often compared
to the complex number $2\pi i$ due to several analogies connecting the multiplicative group $\GG_m$ and the Carlitz's module.
In fact, we have just seen that there are {\em infinitely many} fundamental periods, depending on the various choices of $\vartheta$ we can do
to fix particular maps $\phi_{\CC_\infty,\tau,\vartheta}$.

We can also fix $\vartheta$ and replace $\tau$ by $\tau^n$ for $n>1$. The above discussion applies because the formulas (\ref{dis})
still hold, up to replacement of $q$ by $q^n$. 

\begin{Remark}\label{hellegouarch}
{\em In \cite{Hell1,Hell2} as well as in other references collected in \cite{Hell1}, the reader can find the development of these ideas in yet another framework that
we can consider as the {\em Example (4).} Let us choose, following Hellegouarch, a field $L$ and an
indeterminate $\vartheta$ transcendental over $L$ (including the case $L=\FF_q$). Over the field of formal series $\mathcal{K}=L((1/\vartheta))$,
we consider the $L$-linear endomorphism $\tau$ determined by $$\tau\vartheta=\vartheta^d+a_1\vartheta^{d-1}+a_2\vartheta^{d-2}+\cdots$$ where 
$d$ is an integer $\geq 2$ and $a_1,a_2,\ldots$ are elements of $L$. We then have that $\mathcal{K}^{\tau^n}=L$ for all $n\geq 1$.
Then (see \cite{Hell1}), the series $E$ in the framework of Carlitz module $\phi_{\mathcal{K},\tau,\vartheta}$ induces,
over a suitable field of Puiseux series containing $\mathcal{K}$, a well defined and $L$-linear function,
the kernel of which is a discrete submodule, establishing a good analogue of Theorem \ref{expcarlitz}. Among several
other facts we will see, as a consequence of the results in Section \ref{31}, that the bound $d\geq 2$ is best possible.}
\end{Remark}

\subsection{The Anderson-Thakur function}

Since $\CC_\infty$ is algebraically closed, coherent sequences $\Xi$ in the framework of Example (2) (with a choice of $\vartheta$) exist in $\CC_\infty$
and can be easily constructed by solving iterated Artin-Schreier equations. Just as
for Euler's gamma function (\ref{gammamellin}), there is here a 
canonical choice of coherent sequence which can be made thanks to the exponential function, when $|\vartheta|>1$.

Indeed, theorem \ref{expcarlitz} allows us to immediately verify that if $|\vartheta|>1$,
$$\Xi=\left(\exp_{\tau,\vartheta}\left(\frac{\widetilde{\pi}_{\tau,\vartheta}}{\vartheta^{i+1}}\right)\right)_{i\geq 0}$$
is well defined and is a coherent sequence. We will call the Akhiezer-Baker function $\omega_{\CC_\infty,\tau,\vartheta,\Xi}$ 
the {\em Anderson-Thakur} function associated to the data $(\CC_\infty,\tau,\vartheta)$, assigning to it the simpler notation $\omega_{\CC_\infty,\tau,\vartheta}$. It can be defined, equivalently, by the formulas:
\begin{equation}\label{formulasomega}\omega_{\CC_\infty,\tau,\vartheta}(t)=\exp_{\tau,\vartheta}\left(\frac{\widetilde{\pi}_{\tau,\vartheta}}{\vartheta-t}\right)=\sum_{i\geq 0}\frac{\widetilde{\pi}_{\tau,\vartheta}^{q^i}}{d_i(\vartheta)(t-\vartheta^{q^i})},
\end{equation} where $\exp_{\tau,\vartheta}$ must be considered as an $\FF_q((t))$-linear operator (cf. \cite{Pe}).
From now on, we will also denote by $\omega$ the {\em original Anderson-Thakur function} appearing in \cite{AT}, thus associated with the choice $\vartheta=\theta$.
 We recall that $\omega_{\CC_\infty,\tau,\vartheta}$ satisfies, after
Proposition \ref{PsiX}, the identities:
\begin{equation}\label{eqatomega}
\phi_{\CC_\infty,\tau,\vartheta}(a)\omega_{\CC_\infty,\tau,\vartheta}(t)=a(t)\omega_{\CC_\infty,\tau,\vartheta}(t),\quad a\in\mathcal{A}.\end{equation}
The next proposition contains the fundamental properties of the functions $\omega_{\CC_\infty,\tau,\vartheta}$. They are simple generalizations 
of well known properties of the function $\omega=\omega_{\CC_\infty,\tau,\theta}$ that can be found in \cite{An,AT,ABP,Pe}. The proof of the Proposition
uses the formulas (\ref{formulasomega}) and does not offer additional 
difficulty so it will be omitted.
\begin{Proposition}\label{propomegafundamental}
Let $\vartheta$ be an element of $\CC_\infty$ such that $|\vartheta|>1$. We have the following properties.
\begin{enumerate}
\item The formal series $\omega_{\CC_\infty,\tau,\vartheta}(t)\in\CC_\infty[[t]]$ converges for $t\in\CC_\infty$ such that
$|t|<|\vartheta|$ to a rigid analytic function.
\item The function $1/\omega_{\CC_\infty,\tau,\vartheta}$ extends to an entire rigid analytic function $\CC_\infty\rightarrow\CC_\infty$
with zeros precisely located at the elements $\vartheta,\vartheta^q,\ldots$. 
\item The function $\omega_{\CC_\infty,\tau,\vartheta}$ extends to
a meromorphic function $\CC_\infty\rightarrow\CC_\infty$ with no zeros, and simple poles at $\vartheta,\vartheta^q,\ldots$
of residues $-\widetilde{\pi}_{\tau,\vartheta},-\frac{\widetilde{\pi}_{\tau,\vartheta}^q}{d_1(\vartheta)},\ldots$, where the coefficients $d_i(\vartheta)$ are defined in (\ref{dis}).
\end{enumerate}
\end{Proposition}

\subsection{Functional relations}\label{functionalrelations}

The results of Section \ref{funcrel} can be applied to describe the multiplication relations and the cyclotomic relations for
the function $\omega$. Following \cite{Pe}, we observe that $$\omega(t)=(-\theta)^{\frac{1}{q-1}}\prod_{i>0}\left(1-\frac{t}{\theta^{q^{i}}}\right)^{-1},$$
with a suitable choice of $(q-1)$-th root of $-\theta$.
More generally, if 
$\zeta\in\FF_q^\times$ is such that $\zeta^n=1$ and $\zeta^k\neq 1$ for $0<k<n$, we have
\begin{eqnarray*}
\omega_{\CC_\infty,\tau^n,\theta}(t)&=&(-\theta)^{\frac{1}{q^n-1}}\prod_{i>0}\left(1-\frac{t}{\theta^{q^{ni}}}\right)^{-1},\\
\omega_{\CC_\infty,\tau,\zeta\theta}(t)&=&(-\zeta\theta)^{\frac{1}{q-1}}\prod_{i>0}\left(1-\frac{t}{\zeta\theta^{q^{i}}}\right)^{-1},\\
\omega_{\CC_\infty,\tau,\theta^n}(t)&=&(-\theta^n)^{\frac{1}{q-1}}\prod_{i>0}\left(1-\frac{t}{\zeta\theta^{nq^{i}}}\right)^{-1}.
\end{eqnarray*}
These formulas immediately imply that $\lambda=\mu=1$ in Propositions \ref{propmultiplication} and \ref{propcyclotomic}.
More explicitly, we derive the following Proposition.
\begin{Proposition}
With $\zeta,n$ as above, the following formulas hold.
\begin{enumerate}
\item Multiplication relations. $$\omega_{\CC_\infty,\tau,\theta}(t)=\prod_{i=0}^{n-1}\tau^i\omega_{\CC_\infty,\tau^n,\theta}(t),$$
\item Cyclotomic relations. $$\omega_{\CC_\infty,\tau,\theta^n}(t^n)=\prod_{i=0}^{n-1}\omega_{\CC_\infty,\tau,\zeta^i\theta}(t).$$
\end{enumerate}
\end{Proposition}
The latter relations are very similar to the {\em reflection formulas} for the {\em geometric gamma function} as in Thakur's article \cite{Tha1}. It is for this 
reason that in \cite{Pe}, we have mentioned them as the ``analogues of the reflection formula" for Euler's gamma function.

\subsection{The function $\omega$ and the torsion}

Anderson-Thakur's function $\omega$ interpolates many torsion elements for the Carlitz module in the framework of Example (2), 
just as the gamma function in Proposition \ref{sometorsion}. Indeed, the following Proposition holds.

\begin{Proposition} Let $a$ be irreducible of degree $d$ in $\FF_q[\theta]$. 
Let $\xi_1,\ldots,\xi_d$ be the roots of $a$ in $\FF_{q^d}$. The following identity holds:
$$\omega(\xi_1)+\cdots+\omega(\xi_d)=\exp\left(\frac{\widetilde{\pi}a'}{a}\right),$$ where $a'$ indicates the first derivative of 
$a$ with respect to $\theta$.\end{Proposition}

\noindent\emph{Proof.} We have the following identities:
\begin{eqnarray*}
\sum_{j=1}^d\omega(\xi_j)&=&\sum_{n\geq 0}\exp\left(\frac{\widetilde{\pi}}{\theta^{n+1}}\right)(\xi_1^n+\cdots+\xi_d^n)\\
&=&\sum_{n\geq 0}\exp\left(\frac{\widetilde{\pi}}{\theta^{n+1}}(\xi_1^n+\cdots+\xi_d^n)\right)\\
&=&\exp\left(\widetilde{\pi}\sum_{j=1}^d\sum_{n\geq 0}\theta^{-n-1}\xi_j^n\right)\\
&=&\exp\left(\widetilde{\pi}\sum_{j=1}^d\frac{1}{\theta-\xi_r}\right)\\
&=&\exp\left(\frac{\widetilde{\pi}a'}{a}\right).
\end{eqnarray*}
The essential step is that $\sum_j\xi_j\in\FF_q$ so that
we can use $\FF_q$-linearity of $\exp$.\CVD

\begin{Remark}{\em More generally, it is possible to show the following property, independently noticed by B. Angl\`es. Let $L$ be the compositum 
in $\CC_\infty$ of the field $\FF_q^{\text{alg}}$ and the
various torsion subfields $K(\text{ker}(\phi(\mathfrak{p})))$ for $\mathfrak{p}$ monic irreducible polynomial of $A$. Then, $L$ is equal to the field 
generated over $\FF_q^{\text{alg}}$ by the values $\omega(\xi)$, with $\xi$ varying in $\FF_q^{\text{alg}}$. 
By \cite[Corollary 5]{Pe}, if $\xi$ is in $\FF_{q^d}\setminus\FF_{q^{d-1}}$, then the set of solutions in $\CC_\infty$ of the 
equation $$X^{q^d-1}=(\xi-\theta^{q^d})\cdots(\xi-\theta)$$
is the subset $$\FF_{q^d}^\times\omega(\xi).$$
B. Angl\`es also pointed out that
the function $\omega$ can be understood as an ``universal Gauss-Thakur sum".}\end{Remark}

\section{Back to Example (1): the series $L_x$ as operators}\label{31}

In the specific case of Example (1), we want to investigate the nature of the functions $E,L$.
For this, we need appropriate spaces on which these operators act, but the problem to find natural ones seems to be difficult, unlike the case of Example (2).
Here, we propose an algebra $\mathbb{B}$ of functions over which $L_x$ induces an automorphism for all $x\in\CC$. This is an important 
difference compared with the case of Example (2) with, say, $\vartheta=\theta$, where the exponential function $\exp$ has non-trivial kernel (a free $A$-module of rank one, generated
by Carlitz's period $\widetilde{\pi}$).

\begin{Definition}{\em 
Let us consider the ring $\mathcal{B}$ of functions defined and continuous for $s$ with $\Re(s)\geq 0$,
holomorphic for $\Re(s)>0$,
such that, in any horizontal half-strip $$\{z\in\CC;0<\alpha_0<\Re(z),\beta_0<\Im(z)<\beta_1\},$$
the function $|f(z)|$ is dominated by $c_1e^{c_2|z|}$, where $c_1,c_2$ are positive real numbers depending on the half-strip. 

The ring $\mathcal{B}$ is an integral domain and is an algebra
over the ring of entire functions which are periodic of period $1$. The $F^{\tau}$-algebra $\mathbb{B}=\mathcal{B}\otimes F^\tau$
will be sometimes called the {\em algebra of test functions}.
}\end{Definition}

The algebra $\mathbb{B}$ is endowed with a semi-norm $\|\cdot\|$ in the following way. Let us first define $\|\cdot\|$ over 
the algebra $\mathcal{B}$. Let $f$ be an element of $\mathcal{B}$; we set:
$$\|f\|=\inf\{c>0,\lim_{\alpha\rightarrow\infty}f(\alpha)c^{-\alpha}=0\}.$$
The above, is a semi-norm. It is not a norm; in particular, functions $f$ which are non-constant with $\|f\|=0$
exist, like the function $1/\Gamma(s)$. 
Since for $f,g\in\mathcal{B}$ we have $\|f+g\|\leq\max\{\|f\|,\|g\|\}$, by setting $\|h\|=1$ for all $h\in F^{\tau}\setminus\{0\}$,
there exists a unique semi-norm $\|\cdot\|$ on $\mathbb{B}$ extending the previous one. 
The automorphism $\tau$ of $F$ induces an endomorphism of $\mathbb{B}$. On the other hand,
by the fact that polynomials with complex coefficients belong to $\mathbb{B}$, Carlitz's module $\phi_{F,\tau,s}$ induces an
action of $\mathcal{A}$ over $\mathbb{B}$.

\begin{Remark}{\em 
The algebra $\mathcal{B}$ contains, as a subring, the ring $\mathcal{E}$ of the {\em entire functions of exponential type}
(\footnote{An entire functions of exponential type is an entire 
function $f:\CC\rightarrow\CC$ such that $$|f(s)|\leq c_1e^{c_2|s|},\quad \text{ for all }s\in\CC,$$ for some $c_1,c_2>0$.}).
If $f\in\mathcal{E}$, then $\|f\|$ is less than the exponential of the {\em type} of $f$. We recall that Siegel's $E$-functions are 
examples of entire functions of exponential type. It can also be proved that the function $1/\Gamma(s)$ and 
the function $(1-s)\zeta(s)$ belong to $\mathcal{E}$. But notice that no one of the functions $\Gamma(s+z)$, $z\in\CC$
belongs to $\mathbb{B}$, as Stirling asymptotic estimate indicates.}
\end{Remark}


We now look at $E$ and, more generally, at $L_x$, as operators on $\mathbb{B}$.

\begin{Proposition}\label{autom}
For all $x\in\CC$, $L_x$ defines an $F^\tau$-linear automorphism of 
$F^\tau$-vector spaces
$$L_x:\mathbb{B}\rightarrow\mathbb{B}$$
of inverse $L_{-x}$. 
For all $f\in\mathbb{B}$ and $x\in\CC$, 
$$\|L_x f\|\leq \|f\|.$$
\end{Proposition}
\noindent\emph{Proof.} First of all, for any $x\in\CC$ and $f\in\mathbb{B}$, $L_xf$ is clearly well defined.
Moreover, if $|f(\alpha+i\beta)|\leq c_1e^{c_2\alpha}$ for $\alpha$ big enough, then
$$|L_xf(\alpha+i\beta)|\leq c_1e^{e^{c_2}|x|}e^{c_2\alpha},$$
so that $L_xf\in\mathbb{B}$ and $\|L_xf\|\leq \|f\|$.

The function ${L}_x$ is $F^\tau$-linear because on one side it is additive, 
and on the other side, if $f\in\mathbb{B}$ and $\lambda\in F^\tau$, then:
\begin{eqnarray*}
{L}_x(\lambda f)&=&\sum_{n\geq 0}x^n\frac{\tau^n(\lambda f)}{n!}\\
&=&\lambda\sum_{n\geq 0}x^n\frac{\tau^n(f)}{n!}\\
&=&\lambda{L}_x(f).
\end{eqnarray*}
 
By the fact that ${L}_x{L}_{-x}={L}_{-x}{L}_x=\tau^0$ we see
that if $f$ is such that ${L}_xf=0$, then $f=0$ so ${L}_x$ is bijective.
In particular, ${L}={L}_1$ and $E={L}_{-1}$ induce two 
$F^\tau$-linear automorphisms of $\mathbb{B}$, one inverse of the other. \CVD

\begin{Remark}{\em 
We can also compare the above observations with Remark \ref{hellegouarch} by considering, instead of the field $\mathcal{K}=F$,
the field $\mathcal{K}=\CC((s^{-1}))$ which is also endowed with the automorphism $\tau$ defined by $\tau(s)=s+1$
(so we are in the case $d=1$). The reader can verify that the corresponding functions $L_x$ are all well defined 
$\CC$-automorphisms.}\end{Remark}

\subsection{Compatibility with Mellin's transform}

The use of the algebra $\mathbb{B}$ has some advantages, but, as we pointed out, $\Gamma$ does not count among its elements.
However, the Mellin transform is well defined on $\mathcal{B}$, the operators $L_x$ have good compatibility 
with it, and the gamma function is the Mellin transform of $e^{-t}$.

We focus on the structure of $\CC$-vector space of $\mathcal{B}$, that we filter according to the types.
For $\alpha$ real, let $\mathcal{B}_\alpha$ be the vector space of functions $f$ of $\mathcal{B}$ such that
$|f(t)|\leq c_1e^{-t\alpha}$ for $t$ big enough. We have, for $\alpha\leq \beta$, $\mathcal{B}_{\beta}\subset \mathcal{B}_\alpha$.

If $f(t)$ is a continuous function of the real variable $t\in]0,\infty[$ such that 
$f$ is right continuous 
at $0$ and such that $f(t)=O(e^{-\alpha t})$ for $\alpha>0$ as $t\rightarrow\infty$,
then,
the integral (the Mellin transform of $f$)
$$F(s)=\int_0^\infty t^{s-1}f(t)dt$$
converges in the half-plane $\Re(s)>0$.
For $\alpha>0$, let us denote by $V_\alpha$ the image under the Mellin transform of the space $\mathcal{B}_\alpha$.
We also recall that the Mellin transform, as an operator on a space $\mathcal{B}_\alpha$ for fixed $\alpha$, is injective.

\begin{Lemma}\label{mellin} Let $f$ be an element of $\mathcal{B}_\alpha$ with $\alpha>0$, so that its Mellin transform $F(s)$ defines
an holomorphic function for $\Re(s)>0$. Then, if $x$ is a complex number such that
$$|x|<\alpha,$$
the evaluation of the formal series $L_x\in\QQ[x][[\tau]]$ at $F$ is well defined and is the Mellin 
transform of the function $e^{tx}f(t)\in\mathcal{B}_{\alpha-|x|}$.
\end{Lemma}
\noindent\emph{Proof.} Let $\sigma$ be the real part of $s$. There exists a constant $c>0$ such that, for $n\geq 0$,
$$\left|\frac{x^n}{n!}\int_{0}^\infty t^{s+n-1}f(t)dt\right|\leq c\alpha^{-\sigma}\left|\frac{x}{\alpha}\right|^n\frac{\Gamma(n+\sigma)}{n!}.$$
Under the above conditions and setting $X=|x/\alpha|$, the series 
$$\sum_{n\geq 0}X^n\frac{\Gamma(n+\sigma)}{n!}$$
converges to the function
$\frac{\Gamma(\sigma)}{(1-X)^\sigma}$
and the following transformations hold:
\begin{eqnarray*}
L_x\int_0^\infty t^{s-1}f(t)dt&=&\int_0^\infty L_x(t^{s-1})f(t)dt\\
&=&\int_0^\infty t^{s-1}e^{tx}f(t)dt.
\end{eqnarray*}\CVD

\noindent\emph{Examples.} We have $\Gamma\in V_1$. Therefore, $$L_x\Gamma(s)=\frac{\Gamma(s)}{(1-x)^{s}}\in V_{1-|x|}$$ for all $x$ complex such that
$|x|<1$. In the next table we provide a list of useful identities involving the operators $L_x$.

\medskip

\begin{center}
\begin{tabular}{|l|l||l|}
\hline
$f$ & $L_x(f)$ & condition\\
\hline
\hline
$L_{-x}(f)$ & $f$ &\\
\hline
$t^sf$ & $e^{xt}t^s L_{xt}(f(s))$ & \\
\hline
$sf$ & $sL_x(f)+L_x(\tau f)$ &\\
\hline 
$\tau f$ & $\frac{d}{dx}\left(L_x(f)\right)$ &\\
\hline
\hline 
$\int_{0}^\infty t^{s-1}g(t)dt$ & $\int_{0}^\infty t^{s-1}e^{tx}g(t)dt$ & $g\in\mathcal{B}_\alpha$, $\alpha>|x|$ \\
\hline 
$Q$ & $\sum_{i\geq0}E_i(Q)x^i$ & $Q\in\CC[s]$\\
\hline
$\Gamma(s)$ & $(1-x)^{-s}\Gamma(s)$ & $|x|<1$ \\
\hline
\end{tabular}
\end{center}

\medskip
 
More generally, the operators $L_x$ can be used to construct various classical special functions. 
For example, we can obtain Gauss' hypergeometric functions as images of fractions involving products
of the gamma function via $L_x$.
The reader can find many other examples by just examining any textbook of special functions like \cite{Abra}.

\begin{Remark}{\em We recall a formula which can be deduced from Lommel's Theorem for the Bessel functions (cf. \cite[p. 141]{Bessel}):
$$\sum_{n\geq 0}\frac{(-1)^n(z/2)^n}{n!}J_{\nu-m}(z)=0,\quad \Re(\nu)>0,\quad z\in\CC^\times,$$ where $J_\nu(z)$ denotes 
Bessel's function of first kind. The formula can be rewritten as
$${E}^*((-z/2)^{-s}J_{s}(z))=0,$$ where $E^*=\sum_{n\geq 0}\tau^{-n}/n!$ is the exponential series $E$
associated to the data $(F,\tau^{-1},s)$ and $z$ is fixed. In particular, the kernel of the exponential $E^*$ is non-trivial.}\end{Remark}

\begin{Remark}{\em The author owes the following remark to D. Goss. 
The integral formula of the gamma function (\ref{gammamellin}) can be interpreted as follows.
For $z\in\CC$ such that $\Re(z)>0$, let us consider the family of measures $d\mu_z$ on $\RR_{>0}$ by $d\mu_z=t^z dt/t$. Then 
$$
\Gamma (z)=\int_{\RR_{>0}} e^{-t} d\mu_z(t).
$$

Now, in the framework of Example (2) with $\vartheta=\theta$, let $\delta_a$ be the Dirac measure at $a\in\CC_\infty$ and consider, for a 
parameter $t\in\CC_\infty$, the distribution 
$$
d\nu_t=\sum_{i\geq 0} t^i \delta_{1/\theta^{i+1}}
$$
For $|t|$ small, this converges to a measure on the ring $\mathcal{O}$ of integers of $K_\infty$.
Then, we have
$$
\omega (t)=\int_{\mathcal{O}} \exp(\widetilde{\pi}z) d\nu_t(z),
$$ similar to (\ref{gammamellin}). Moreover, for $a\in A$,
\begin{eqnarray*}
a(t)\omega (t)&=&\phi(a)\int_{\mathcal{O}} \exp(\widetilde{\pi}z) d\nu_t(z)\\
&=&\int_{\mathcal{O}} \phi(a)\exp(\widetilde{\pi}z) d\nu_t(z)\\
&=&\int_{\mathcal{O}} \exp(\widetilde{\pi}az) d\nu_t(z),
\end{eqnarray*}
in analogy with the formula
$$\nu^{-s}\Gamma(s)=\int_{\RR_{>0}} e^{-\nu t} d\mu_z(t),\quad \Re(\nu)>0.$$

However, the analogy is only partial. In (\ref{gammamellin}),
the function we integrate is related to the exponential function $e^x=\sum_{i\geq0}x^i/i!$ of the multiplicative group $\GG_m$
but not to the exponential function of the generalized Carlitz module of Example (1). In fact, the function
$e^x$ cannot be naturally seen as the exponential function of a generalized Carlitz's module as described in this text, for the
reason that $\ZZ$ (the ring that here should correspond to $\mathcal{A}$) is not of the form $L[\vartheta]$ for a field $L$ and 
an element $\vartheta$ transcendental over $L$!}\end{Remark}

\section{Some applications of Carlitz formalism}

We return here to the framework of Example (1) and observe that some simple and classical properties of Hurwitz's zeta function
$$\zeta(z,s)=\sum_{n\geq 0}\frac{1}{(n+s)^{z}},\quad \Re(z)>1,\quad \Re(s)>0$$
can be easily verified 
with the Carlitz formalism and some standard asymptotic estimates.

Riemann's integral formula
\begin{equation}\label{riemann}
\zeta(z)\Gamma(z)=\int_0^\infty t^{z-1}\frac{1}{e^t-1}dt\end{equation}
extends to Hurwitz zeta function 
by means of the integral formula:
\begin{equation}\label{hurwitz}
\zeta(z,s)\Gamma(z)=\int_0^\infty e^{-s t}t^{z-1}\frac{1}{1-e^{-t}}dt
\end{equation}
(Riemann's formula corresponds to the value $s=1$).
This implies that $\zeta(z,s)$, for fixed $s$, extends to a meromorphic 
function on $\CC$ with, as unique singularity, a simple pole at $z=1$ with residue $1$ (thus independent
on the choice of $s$).

The next question is then to look at the Laurent series expansion of $\zeta(z,s)$ at $z=1$. The coefficients are 
often used to define the so-called Stieltjes constants. Here is a very classical result for the constant term
in this Laurent series.
\begin{Proposition}\label{proplimitatone} The following limit holds:
$$\lim_{z\rightarrow 1}\zeta(z,s)-\frac{1}{z-1}=-\frac{\Gamma'(s)}{\Gamma(s)}.$$
\end{Proposition}
\noindent\emph{Proof.} This limit can be computed in a variety of ways. We propose here to use Carlitz's formalism. We begin by showing that
the function
$$f(s):=\Gamma(s)\left(\lim_{z\rightarrow 1}\zeta(z,s)-\frac{1}{z-1}\right)$$
is an element of $s^2$-torsion.

We compute first $\phi(s)f(s)$, that is, $sf(s)-f(s+1)$. We find:
\begin{eqnarray*}
\phi(s)f(s)&=&s\left(\Gamma(s)\left(\lim_{z\rightarrow 1}\zeta(z,s)-\frac{1}{z-1}\right)\right)-
\Gamma(s+1)\left(\lim_{z\rightarrow 1}\zeta(z,s+1)-\frac{1}{z-1}\right)\\
&=&s\Gamma(s)\left(\lim_{z\rightarrow 1}\zeta(z,s)-\zeta(z,s+1)\right)\\
&=&\Gamma(s).
\end{eqnarray*}
Therefore, $\phi(s^2)f(s)=\phi(s)\Gamma(s)=0$, so $f(s)$ is an element of $s^2$-torsion as claimed.
This means, by Proposition \ref{torsion}, that there exist two functions $a(s),b(s)\in F^\tau$ such that
\begin{equation}\label{ab}
f(s)=a(s)\Gamma(s)+b(s)\Gamma'(s),
\end{equation}
and we need show that $a,b$ are constant, with $a=0$, $b=-1$. 
For this, an asymptotic estimate will suffice. We appeal to Euler-MacLaurin summation formula which yields, for $\Re(z)>1$ and $\Re(s)>0$:
$$\zeta(z,s)-\frac{1}{z-1}=\frac{s^{1-z}-1}{z-1}+\frac{s^{-z}}{2}-z\int_0^\infty\frac{t-[t]-\frac{1}{2}}{(t+s)^{z+1}}dt.$$
Now, the limit for $z\rightarrow 1$ of $\frac{s^{1-z}-1}{z-1}$ is well defined, being equal to $-\log(s)$ so that we have the 
well defined limit
$$\lim_{z\rightarrow1}\zeta(z,s)-\frac{1}{z-1}=-\log(s)+\frac{1}{2s}-\int_0^\infty\frac{t-[t]-\frac{1}{2}}{(t+s)^{2}}dt.$$
For $s$ real tending to $\infty$, the left-hand side of the above identity is then asymptotically equivalent to 
$-\log(s)$. But it is easy to verify that the only element of the $F^\tau$-span of the functions
$1,\frac{\Gamma'(s)}{\Gamma(s)}$ with this property is precisely 
$-\Gamma'(s)/\Gamma(s)$ so that $a=0$ and $b=-1$.\CVD

\begin{Remark}{\em More generally, 
let $\Psi(s,t)$ be the function
$$\Psi(s,t)=\sum_{n\geq 1}\zeta(n+1,s)t^n.$$
Then, the identity
$$\Psi(s,t)=\psi(s)-\psi(s-t)$$
holds, for $\Re(s)>\Re(t)>0$. The proof uses the fact that $\psi(s)-\Psi(s,t)$ is solution of the non-homogeneous $\tau$-difference equation
(\ref{taudigamma}), and standard asymptotical estimates. This yields well known formulas (see for instance \cite[25.11.12]{Abra}):
\begin{equation}\label{propoidentity2}\zeta(n+1,s)=\frac{(-1)^{n+1}\psi^{(n)}(s)}{n!},\quad n\geq 1.
\end{equation}}\end{Remark}

\begin{Remark}{\em It is also well known that the above results can be directly deduced from Weierstrass' product expansion of the gamma function.
Carlitz's formalism and asymptotic estimates allow to deduce this product expansion as well; this is implicit in the enlightening introduction of Aomoto and Kita's book
\cite[p. 1-3]{AKI}.}\end{Remark}

We provide a last application of Carlitz formalism, this time involving the operator $L_x$.
The following formula appears in \cite{Wilton} and in \cite[p. 137 (6.6)]{Srivastava}, and generalizes a classical formula by Landau. 
If $x,\alpha$ are complex numbers such that $|x|<|\alpha|$ and $\alpha$ is not in the set $\{0,-1,-2,\ldots\}$,
then, for $s$ complex number different from $1$,
$$\sum_{n\geq 0}\binom{s+n-1}{n}\zeta(s+n,\alpha)x^n=\zeta(s,\alpha-x).$$ 
In our notations, we have:
\begin{Proposition}\label{srivastava}
For $x,\alpha$ complex such that $|x|<|\alpha|$ and $\Re(\alpha)>0$, we have:
$${L}_x(\Gamma(s)\zeta(s,\alpha))=\Gamma(s)\zeta(s,\alpha-x).$$
\end{Proposition}
\noindent\emph{Proof.} It can be deduced from (\ref{hurwitz}) and an appropriate variant of Lemma \ref{mellin} when $\alpha>0$ but it extends 
to complex $\alpha$ with $\Re(\alpha)>0$ as well. The following intermediate identities can be easily verified:
\begin{eqnarray*}
{L}_x(\Gamma(s)\zeta(s,\alpha))&=&\int_0^\infty e^{-\alpha t}{L}_x(t^{s-1})\frac{1}{1-e^{-t}}dt\\
&=&\int_0^\infty e^{-\alpha t}t^{s-1}\frac{e^{tx}}{1-e^{-t}}dt\\
&=&\int_{0}^\infty e^{-(\alpha-x)t}t^{s-1}\frac{1}{1-e^{-t}}dt\\
&=&\zeta(s,\alpha-x)\Gamma(s).
\end{eqnarray*} \CVD

\section{Link with the functions $L(\chi_t^\beta,\alpha)$.}

We end our discussion about Carlitz formalism by considering once again Example (2) and we restrict our attention to the case $\vartheta=\theta$. 
Let $A^+$ be the set of monic polynomials of $A=\FF_q[\theta]$. In \cite{Pe}, we have introduced, for 
$\beta\geq 0$ and $\alpha>0$, the series:
$$L(\chi_t^\beta,\alpha)(t)=\sum_{a\in A^+}\chi_t(a^\beta)a^{-\alpha}\in K_\infty[[t]].$$ 
Here, for a given polynomial $a$ of $A$, $\chi_t(a)$ denotes the polynomial function obtained 
by substituting $\theta$ by $t$. In fact, the series $L(\chi_t^\beta,\alpha)$ converges for $t\in\CC_\infty$ with $|t|\leq 1$ and
has analytic extension to the whole $\CC_\infty$-plane.

We recall, following \cite{Pe}, that for $t\in\{\theta^{q^k},k\in\ZZ\}$, $L(\chi_t^\beta,\alpha)(t)$ specializes to 
a {\em Carlitz zeta value}. For example, for $t=\theta$, we get
$$L(\chi_t^\beta,\alpha)(\theta)=\zeta(\alpha-\beta)=\sum_{d\geq 0}\sum_{a\in A^+(d)}a^{\beta-\alpha},$$ where 
$A^+(d)$ denotes the set of monic polynomials of degree $d$.

Also, in \cite{Pe}, we have stressed that special values of $L$-series associated to a certain Dirichlet character 
may also be obtained as special values of the series $L(\chi_t^\beta,\alpha)$, this time substituting $t$ with a root of unity. This 
already gives the series $L(\chi_t^\beta,\alpha)$ a flavor of 
Hurwitz's zeta function because the latter function is also used to
compute $L$-series associated to Dirichlet characters in the domain of their analytic 
continuation.

There are other signs which allow to pursue the analogy with the Hurwitz's zeta function. One of the results in \cite{Pe} is the following formula:
\begin{equation}
L(\chi_t,1)(t)=-\frac{\widetilde{\pi}}{(t-\theta)\omega(t)}.
\end{equation}
Among several things, this formula, together with Proposition \ref{propomegafundamental}, allows
to prove directly that $L(\chi_t,1)(t)$ can be analytically extended to $\CC_\infty$, and provides
important information on certain Carlitz's zeta values, as well as on certain trivial zeros of 
Goss' zeta function (see \cite{Pe}). Now, this goes in the same direction of 
Proposition \ref{proplimitatone} and (\ref{propoidentity2}).

Here is a last analogy. In a work in progress, the author and R. Perkins recently 
considered an extension of the semi-character $\chi_t$ to the whole $\CC_\infty$-plane.
In order to construct this, one appeals to Anderson-Thakur's function by setting, for 
$z\in\CC_\infty$,
$$\chi_t(z)=\frac{\phi(z)\omega(t)}{\omega(t)},$$ $\phi(z)$ being the operator of $\CC_\infty[[\tau]]$ 
as in Definition \ref{definitioncarlitz}. The evaluation at $\omega$ is possible (there is a suitable notion of radius of convergence
to take into account here). In particular, if $z=a\in A$,
we get the above semi-character by (\ref{eqatomega}). Now, extending $\tau$ to the 
ring $\CC_\infty[[t,t_1]]$ with the rule $\tau t_1=t_1$, where $t_1$ is another variable, we can
as well evaluate the operator $\log=\sum_{n\geq 0}(-1)^nl_n^{-1}\tau^n$ at a series 
of $\CC_\infty[[t,t_1]]$. We show that \begin{equation}
\log\left(\omega\sum_{a\in A^+}\chi_t(a^{-1})\chi_{t_1}(a)\right)=\frac{-\widetilde{\pi}}{t-\theta}L(\chi_{t_1},1),\label{srivastavaform}
\end{equation}
a formula which has some similarity with Proposition \ref{srivastava}.

\begin{Remark}{\em The theory of the generalized Carlitz's module can be extended to higher rank case in parallel with the theory of Drinfeld modules. 
We will not discuss it here, but 
again the classical theory of special functions offers abundant examples of these issues. For example, the well known {\em contiguity relations} for 
Gauss' hypergeometric functions may be viewed, after suitable normalizations as particular cases of {\em generalized Drinfeld modules of rank two} in the framework of Example (1).
The author has already partially tracked this analogy in the recent paper \cite{archiv}.}\end{Remark}

\medskip

\noindent{\bf Acknowledgements.} The author is indebted to D. Goss for enriching discussions and encouragement in writing this work, and  
to B. Angl\`es, D. Essouabri and D. Goss for several useful comments on earlier 
versions of this paper. The author is also thankful to the referee for very helpful commentaries that allowed to improve this text.


\begin{thebibliography}{9}

 \bibitem{Abra} M. Abramovitz \& I. A. Stegun. {\em Handbook of Mathematical Functions.}
 
 \bibitem{Akh} N. I. Akhiezer. {\em Elements of the theory of elliptic functions.} Izdat.
Nauka, Moscow 1970. English transl., American Mathematical Society, 1990.

\bibitem{An} G. Anderson. {\em $t$-motives}, Duke Math. J. 53 (1986), 457-502. 



\bibitem{AT} G. Anderson \& D. Thakur. {\em Tensor powers of the Carlitz module and zeta values.} Ann. of Math. 132 (1990), 159-191.


\bibitem{ABP} G. Anderson, D. Brownawell \& M. Papanikolas, {\em Determination of 
the algebraic relations among special $\Gamma$-values in positive characteristic,} Ann. of Math. 160 
(2004), 237-313. 


\bibitem{AKI} K. Aomoto \& M. Kita. {\em Theory of Hypergeometric Functions.} Springer Monographs in Mathematics, (2011).


\bibitem{Bak} H. F. Baker. {\em Note on the foregoing paper ``commutative ordinary differential operators" by J. L. Burchnall and T. W.
Chaundy.} Proc. Royal Soc. London (A) 118 (1928). 584-593.

\bibitem{Ca0} L. Carlitz. {\em On certain functions connected with polynomials in a Galois
field.} Duke Math. J. 1, 137Ð168, (1935).

\bibitem{Ca1} L. Carlitz. {\em An analogue of the von Staudt-Clausen theorem.} Duke Math. J. 3, 503-517, (1937).

\bibitem{Ca2} L. Carlitz. {\em An analogue of the Staudt-Clausen theorem.} Duke Math. J. 7, 62-67, (1940).

\bibitem{Ca3} L. Carlitz. {\em An analogue of the Bernoulli polynomials.} Duke Math. J. 8, 405-412, (1941).


\bibitem{ChaHru} Z. Chatzidakis \& E. Hrushovski. {\em Model theory of difference fields}, Trans. of the AMS, 351, pp. 2997-3071, (1999).

\bibitem{Cohn} R. M. Cohn. {\em Inversive difference fields.} Bull. Amer. Math. Soc. Volume 55, Number 6 (1949), 597-603. 
















\bibitem{Go} D. Goss. {\em Basic structures of function field arithmetic.} Ergebnisse der Mathematik und ihrer Grenzgebiete, 35. Springer-Verlag, Berlin, (1996).



\bibitem{Hell1} Y. Hellegouarch. {\em Galois Calculus and Carlitz Exponentials}. In: The Arithmetic of Function Fields, (eds: D. Goss et al), de Gruyter (1992), 33-50.

\bibitem{Hell2} Y. Hellegouarch \& F. Recher. {\em Generalized $t$-modules}. Journal of Algebra 187, 323-372 (1997).






\bibitem{Kri0} I. M. Krichever. {\em Generalized elliptic genera and Baker-Akhiezer functions.} Math. Notes,  Volume 47, Number 2 (1990), 132-142, DOI: 10.1007/BF01156822. Translated from Matematicheskie Zametki, Vol. 47, No. 2, pp. 34Ð45, February, 1990.

\bibitem{Kri} I. M. Krichever. {\em Methods of algebraic geometry in the theory of non linear equations.} Russ. Math. Surv. 32, (1977), 185-213

\bibitem{Levin} A. Levin. {\em Difference Algebra.} Algebra and applications, Vol. 8. Springer, 2008.

\bibitem{Mum} D. Mumford. {\em An algebro-geometric construction of
commuting operators and solutions to the Toda
lattice equation, KdV equation and related nonlinear
equations.} International Symposium on Algebraic
Geometry (Kyoto, 1977) (M. Nagata, ed.), Kinokuniya,
Tokyo, 1978, pp. 115-153.






\bibitem{archiv} F. Pellarin. {\em Estimating the order of vanishing at infinity of Drinfeld quasi-modular forms.}
To appear in J. Reine Angew. Math. published online 02-06-2012.





\bibitem{Pe} F. Pellarin. {\em Values of certain $L$-series in positive characteristic.} Annals of Mathematics 176 (2012), 1-39.

\bibitem{Praagman} C. Praagman. {\em Fundamental solutions for meromorphic linear difference equations in the complex plane, and related problems.} Journal fŸr die reine und anger. Math. 369, 101-109, (1986).




\bibitem{Srivastava} H. M. Srivastava.
{\em A unified presentation of certain classes of series of the Riemann zeta function.} Riv. Mat. Univ. Parma (4), 14, 1988. 1-23.

\bibitem{Tha1} D. Thakur. {\em Gamma functions for function fields and Drinfel'd modules.} Ann. of Math., 2nd Ser., Vol. 134, No. 1. (Jul., 1991), pp. 25-64.


\bibitem{Tha} D. Thakur. {\em Function field arithmetic}. World Scientific, (2004).


\bibitem{Bessel} G. N. Watson. {\em A treatise on the theory of Bessel functions}. Cambridge University press, 1922.

\bibitem{Wilton} 
J . R. Wilton. {\em A proof of Burnside's formula for $\log\Gamma(x + 1)$ and certain allied properties of Riemann's $\zeta$-function.} Messenger Math. 52 (1922/1923), 90-93.

\end{thebibliography}
\end{document}